\documentclass[11pt]{amsart}
\usepackage{amsmath}
\usepackage{amsfonts}
\usepackage{amssymb}
\usepackage{amsthm}
\usepackage{amsrefs}
\usepackage{graphicx}
\usepackage{pb-diagram}
\usepackage{epstopdf}
\usepackage{amscd}
\usepackage{pdfpages}
\usepackage{bbm}
\usepackage{multirow}
\usepackage[mathscr]{eucal}
\usepackage{epsfig,epsf,epic}
\usepackage{pdfsync}
\usepackage{enumerate}
\usepackage{etex}
\usepackage[all]{xy}
\usepackage{fancyref}
\usepackage{mathtools}
\usepackage{scalerel,stackengine}
\usepackage{comment}
\usepackage{hyperref}

\stackMath

\newtheorem{theorem}{Theorem}[section]
\newtheorem{prop}[theorem]{Proposition}

\newtheorem{corollary}[theorem]{Corollary}
\newtheorem{definition}[theorem]{Definition}
\newtheorem{example}[theorem]{Example}

\newtheorem{lemma}[theorem]{Lemma}
\newtheorem{remark}[theorem]{Remark}

\newtheorem{notation}[theorem]{Notations}
\newtheorem{lem-defn}[theorem]{Lemma/Definition}
\numberwithin{equation}{section}

\DeclareMathOperator{\Sym}{Sym}

\newcommand{\A}{\mathcal A}

\newcommand{\W}{\mathcal W}
\newcommand{\C}{\mathbb{C}}

\newcommand{\gqrep}{\beta}

\begin{document}
\title[Symmetry in Deformation quantization and Geometric quantization]{Symmetry in Deformation quantization and Geometric quantization}	
	
\author[Leung]{Naichung Conan Leung}
\address{The Institute of Mathematical Sciences and Department of Mathematics\\ The Chinese University of Hong Kong\\ Shatin \\ Hong Kong}
\email{leung@math.cuhk.edu.hk}
	
\author[Li]{Qin Li}
\address{Institute for Quantum Science and Engineering\\
	Southern University of Science and Technology, Shenzhen, China}
\email{liqin@sustech.edu.cn}

\author[Ma]{Ziming Nikolas Ma}
\address{Department of Mathematics, \\
	Southern University of Science and Technology, Shenzhen, China}
\email{mazm@sustech.edu.cn}

\begin{abstract}
In this paper, we explore the quantization of Kähler manifolds, focusing on the relationship between deformation quantization and geometric quantization. We provide a classification of degree 1 formal quantizable functions in the Berezin-Toeplitz deformation quantization, establishing that these formal functions are of the form \( f = f_0 - \frac{\hbar}{4\pi}(\Delta f_0 + c) \) for a certain smooth (non-formal) function \(f_0\). If \(f_0\) is real-valued then $f_0$ corresponds to a Hamiltonian Killing vector field. In the presence of Hamiltonian \(G\)-symmetry, we address the compatibility between the infinitesimal symmetry for deformation quantization via quantum moment map and infinitesimal symmetry on geometric quantization acting on Hilbert spaces of holomorphic sections via Berezin-Toeplitz quantization. 
\end{abstract}

\maketitle

\section{Introduction}\label{sec:introduction}

Kähler manifolds $(X,\omega, J)$ are central to the study of complex geometry and symplectic geometry due to their unique structures, which encompass both complex and symplectic aspects. These manifolds naturally facilitate the construction of quantum theories, particularly in the context of deformation quantization. In \cite{CLL}, the influence of the Kähler condition on quantum theories is examined by reformulating Kapranov's construction of an $L_\infty$-algebra structure on Kähler manifolds as a flat connection, denoted $D_K$, on the holomorphic Weyl bundle $\W_X$. Futhermore, for every choice of closed $2$-form $\alpha \in  \A_X^{1,1}[[\hbar]]$ of type $(1,1)$, $D_K$ can be quantized to yield a Fedosov abelian connection $D_{\alpha}$, expressed as 
\begin{equation*}
D_{\alpha} = \nabla - \delta + \frac{1}{\hbar}[I_{\alpha}, -]_{\star},
\end{equation*}
 where $\nabla$ represents the Levi-Civita connection, $I_{\alpha}$ is a 1-form valued section of $\W_{X,\C}$. This connection satisfies the Fedosov equation $D_{\alpha}^2=0$, as in the framework of Fedosov's construction of deformation quantization \cite{Fed}. This flat connection is essential for understanding the relationship between geometric structures and quantum properties. For instance, by identifying the sheaf of formal smooth functions $\mathcal{C}^{\infty}(X)[[\hbar]]$ with the sheaf of flat sections under $D_{\alpha}$, mapping a smooth function $f$ to a flat section $O_f$, we obtain a star product $\star_{\alpha}$ induced from the fiberwise star product in $\W_{X,\C}$, deforming the commutative product structure. 

On the other hand, geometric quantization focuses on the Hilbert space $\mathcal{H}_k := H^0_{\bar{\partial}}(X,L^{\otimes k})$ of holomorphic sections of tensor powers of the prequantum line bundle $L$ on $X$ (with $[\omega]=c_1(L)$). It is natural to consider how the deformed non-commutative ring of smooth functions can act on $\mathcal{H}_k$. A prominent example is through the Berezin-Toeplitz operators:
\begin{align*}
	T : C^{\infty }\left( X\right) &\rightarrow End\left( H^0_{\bar{\partial}}(X,L^{\otimes k}) \right)  \\
	f &\mapsto T_{f}=\Pi \circ m_{f},
\end{align*}%
which is given by the multiplication $m_f$ followed by the orthogonal projection $\Pi :\Gamma _{L^{2}}\left( X,L^{\otimes k}\right) \rightarrow H^0_{\bar{\partial}}(X,L^{\otimes k})$, leading to the Berezin-Toeplitz star product $\star_{\text{BT}}$. However, this does not constitute an honest action.

When $X$ and the Karabegov form $K = -\frac{1}{\hbar}\cdot\omega+\alpha$ are real analytic, a sheaf of Bargmann-Fock modules over the Weyl bundle $\W_{X,\C}$ was constructed in \cite{CLL3}, resulting in the Bargmann-Fock sheaf $\mathcal{F}_{\alpha}^{\text{flat}}$, which is naturally equipped with a Fedosov type flat connection $D_{\alpha}$ that is compatible with the module structure. This induces a natural action of the sheaf of analytic functions $\left(C^{\omega}_X[[\hbar]],\star_\alpha\right)$ on the space of flat section of $\mathcal{F}_{\alpha}^{\text{flat}}$ under $D_{\alpha}$, which can be identified with a formal Toeplitz-type action. In \cite{CLL23}, $\alpha$ is further assumed to be a polynomial in $\hbar$, a notion of formal quantizable functions is defined as a dense subsheaf of $C^{\omega}_X[[\hbar]]$ consisting of those $f$'s satisfying a finiteness condition in antiholomorphic part and $\hbar$ in the expression of $O_f$. This makes it possible to evaluate $O_f$ at $\hbar = \frac{1}{k}$ for large enough $k$ and gives a notion of quantizable function of level $k$, this sheaf is denoted by $\mathcal{C}_{\alpha,k}$. Moreover, we can identify flat sections of $\mathcal{F}_{\alpha}^{\text{flat}}$ with holomorphic sections of $L^{\otimes k}$ such that the action of quantizable functions of level $k$ on $L^{\otimes k}$ is further identified with the action of the sheaf of twisted differential operators on $L^{\otimes k}$. In other words, sections of $\mathcal{C}_{\alpha,k}^\infty$ are quantized as twisted differential operators on $L^{\otimes k}$ through this process. 

The simplest class of quantizable functions are those of degree 1 (see Definition \ref{definition:degree_in_weyl}). Our first result is a classification of these quantizable functions.

\begin{theorem}\label{theorem:reformulated_theorem_1}
	For the Berezin-Toeplitz deformation quantization whose Karabegov form is given by \(K_{\text{BT}} = -\frac{1}{\hbar} \cdot \omega +  \text{Ric}_X\), a degree 1 formal quantizable function \(f\) is exactly of the form 
	$$
	f = f_0 - \frac{\hbar}{4\pi}(\Delta f_0 + c),
	$$
	where \(c\) is a constant and \(f_0\) is a smooth function such that the $(1,0)$-component \(V_{f_0}^{1,0}\) of its Hamiltonian vector field $V_{f_0}$ is holomorphic. Evaluating at \(\hbar = \frac{1}{k}\) yields all degree 1 quantizable functions at level \(k\) for sufficiently large \(k\).
\end{theorem}

It turns out that when \(f_0\) is a real-valued function, the condition that \(V_{f_0}^{1,0}\) being holomorphic corresponds exactly to \(V_{f_0}\) being a Killing vector field. Such $f_0$'s form a Lie algebra describing classical infinitesimal symmetry of the K\"{a}hler manifold $(X,\omega,J)$. It turns out that this classical infinitesimal symmetry coincides with that for Berezin-Toeplitz star product $\star_{\text{BT}}$ by the work of \cite{Muller-Neumaier-2}. In other words, quantum infinitesimal symmetry for deformation quantization arises from classical symmetry. 

Classical infinitesimal symmetry of symplectic manifold $(X,\omega)$ is described by the Lie algebra of symplectic vector fields. It is common to further restrict ourselves to Hamiltonian vector field $V = V_f$ obtained from some smooth function $f$, satisfying $V_f(g) = \{f,g\}$ for any smooth function $g$. On the quantum level, we say that a formal vector field $V\in\Gamma(X, TX_{\C})[[\hbar]]$ has a quantum Hamiltonian $f \in \mathcal{C}^{\infty}(X)[[\hbar]]$ associated to it if the equation \(V(g) = \frac{1}{\hbar} [ f, g ]_{\star}\) is satisfied for all $g \in \mathcal{C}^{\infty}(X)[[\hbar]]$, where \([\cdot,\cdot ]_{\star}\) is the commutator of the \(\star_{\alpha}\)-product. We give a characterization of the existence of quantum Hamiltonian $f$ for $V$ for general $\alpha$. 

\begin{theorem}\label{theorem:general_quantum_hamiltonian_reformulated}
	If $V$ is a symplectic vector field preserving $\alpha$, then $f$ is a quantum Hamiltonian associated $V$ if and only if the equation
	$$
		-\frac{2\pi\hbar}{\sqrt{-1}}\cdot\iota_VK=df
	$$
	is satisfied.
	\end{theorem}
In particular for Berezin-Toeplitz quantization, we see that the notion of degree $1$ formal quantizable function $f$ is equivalent to it being a quantum Hamiltonian associated to $V_{f_0}$ when $f_0$ is a real-valued function. This gives a geometric understanding of real-valued degree $1$ formal quantizable function. 
\begin{corollary}\label{theorem:reformulated_corollary_1}
	When \(f \in C^{\infty}(X)[[\hbar]] \) is real-valued, it is a degree $1$ formal quantizable function if and only if it is a quantum Hamiltonian associated to a Killing vector field for Berezin-Toeplitz quantization.
\end{corollary}

Finally, when \((X, \omega, J)\) has a Hamiltonian \(G\)-symmetry with moment map \(\mu : \mathfrak{g}\rightarrow \mathcal{C}^{\infty}(X)\) describing its classical infinitesimal symmetry, this symmetry can be lifted to the Berezin-Toeplitz deformation quantization using the explicit quantum moment map \(\mu_\hbar : \mathfrak{g} \rightarrow C^\infty(X)[[\hbar]]\) defined by the equation
$$
\mu_\hbar(\xi) := \mu(\xi) - \frac{\hbar}{4\pi}\Delta\mu(\xi)
$$ 
for \(\xi \in \mathfrak{g}\) using the above Corollary \ref{theorem:reformulated_corollary_1}. Evaluating at \(\hbar = \frac{1}{k}\) yields \(\mu_k(\xi)\). On the other hand, Hamiltonian $G$-action preserving prequantum line bundle $L$ forms a natural \(G\)-representation, \(G \rightarrow GL(\mathcal{H}_k)\), where the geometric quantization is given by \(\mathcal{H}_k:=H^{0}_{\bar{\partial}}(X,L^{\otimes k})\). This representation describes symmetry for the geometric quantization $\{\mathcal{H}_k\}_{k\in \mathbb{N}}$, which is closely related to the Guillemin-Sternberg conjecture \cite{Guillemin-Sternberg} on the commutativity of symmetry with geometric quantization ( For the proof of the conjecture and other generalizations, readers may see e.g., \cites{Mathai_Zhang, Ma-Zhang, Tian-Zhang, Vergne}). The infinitesimal symmetry is described by the Lie algebra representation 
\begin{equation}\label{equation: G-representation-Hilbert-space}
	\gqrep_k : \mathfrak{g} \rightarrow gl(\mathcal{H}_k).
\end{equation}
Therefore, it is natural to ask whether the infinitesimal symmetry of deformation quantization and that of geometric quantization are compatible. This leads us to consider the following diagram:

$$
\xymatrix{ & C^\infty_{\alpha,k}(X) \ar[dr]^{\mathcal{D}}  & \\
	\mathfrak{g} \ar[ur]^{k\cdot\mu_k}\hspace{2mm}  \ar[rr]_{\beta_k} & & gl(\mathcal{H}_k)
}
$$
where \(\mathcal{D}\) refers to the Bargmann-Fock action. Our last theorem addresses this question.

\begin{theorem}
	The above diagram commutes for Berezin-Toeplitz quantization, i.e. we have the equation $\mathcal{D} \circ \mu_k = \beta_k$. 
\end{theorem}

%As a quantization scheme on general prequantizable K\"ahler manifolds, the Toeplitz operators induce the Berezin-Toeplitz star products: Let $\hslash = \frac{1}{ m}$ and consider the $m$-th power $L^m$, with corresponding Hilbert space $\mathcal{H}^{(m)}=H^0_{L^2}(X,L^{\otimes m})$ and Toeplitz operator $T^{(m)}$. There is the asymptotic expansion as $m \rightarrow \infty$: 
%$$T^{(m)}_f \circ T^{(m)}_g \sim T^{(m)}_{f \star_{\hslash} g},$$
%which uniquely determines $\star_{\hslash}$. The Karabegov form defined in \cite{karabegov2000identification}, classifying star products, is given by $-\frac{1}{\hslash} \omega+ iF_{K}$ for $\star_{\hslash}$. We may consider the quantum moment map $\tilde{\moment}$ for the Karabegov form $\frac{1}{\hslash} \omega - iF_{K}$ as a modification of $\moment$ in the presence of Hamiltonian $G$-symmetry. In general, the above diagram only commute in the asymptotic sense as $\hslash \rightarrow \infty$. See Section \ref{subsection: quantum-moment-map} and Remark \ref{rem:asymptotic_communte} for a more detailed explanation. 

The organization of this paper is as follows: in section \ref{section: Fedosov-quantization}, we will recall the Fedosov quantization method, including both deformation quantization and geometric quantization, and the notion of quantizable functions which can be quantized to holomorphic differential operators on the Hilbert spaces in geometric quantization; in Section \ref{section: first-order-quantizable-function}, we will give a full description of first order quantizable functions. In section \ref{section: quantum-Hamiltonian}, we give a class of examples of first order quantizable functions which arise from symmetries preserving the star product. As an application, we will show in section \ref{section: quantum-moment-map} that for the Berezin-Toeplitz quantization, the quantum symmetry on algebra of observable and Hilbert spaces (geometric quantization) are naturally compatible. 

\subsection*{Conventions}
\begin{itemize}
	\item  Let $X$ be a smooth manifold. We denote by $\Omega^k_X$ the bundle of differential $k$-forms on $X$ and by $\Omega^\bullet_X = \bigoplus_{k}\Omega^k_X$ the full de Rham complex. Global smooth differential forms on  $X$ will be denoted by 
	$$
	\A_X^\bullet=\Gamma(X, \Omega^\bullet_X), \quad \text{where}\ \A_X^k=\Gamma(X, \Omega^k_X).
	$$
	Given a vector bundle $E$, the complex of $E$-valued differential forms is denoted by
	$$
	\A^\bullet_X(E)=\Gamma(X, \Omega^\bullet_X\otimes E). 
	$$
	\item For a complex manifold $X$, we let
	\begin{itemize}
		\item $TX$ and $T^*X$ denote the holomorphic tangent and cotangent bundles respectively;
		\item $\overline{TX}$ and $\overline{T^*X}$ denote the anti-holomorphic tangent and cotangent bundles respectively;
		\item $TX_{\mathbb{R}}$ and $T^*X_{\mathbb{R}}$ denote the real tangent and cotangent bundles respectively;
		\item $TX_\C$ and $T^*X_\C$ denote the complexified tangent and cotangent bundles respectively.
		\item $L$ denotes the pre-quantum line bundle on $X$
		\item $\nabla$ denotes the Levi-Civita connection on $X$, and also its associated connection on the Weyl bundles by abuse of notations. $\nabla^{L^{\otimes k}}$ denotes the connection on $L^{\otimes k}$. 
		\item $\Delta$ denotes the Laplacian on $X$. 
		\item We let $d=\partial+\bar{\partial}$ denote the type decomposition of the de Rham differential on $X$. 
		\item We let $J$ denote the complex structure on $X$. 
		\item For every vector field $V$ on $X$, we let $\mathcal{L}_V$ denote the associated Lie derivative and $\iota_V$ the contraction between vector field with differential forms. 
	\end{itemize}
	\medskip
	
	\item $\mathcal{H}_m$ denotes the space of holomorphic sections of $L^{\otimes m}$.
	
	\item $G$ denotes a compact Lie group acting on $X$, and $\mathfrak{g}$ denotes the Lie algebra of $G$. For any $\xi\in\mathfrak{g}$, we let $V_\xi$ denote the vector field on $X$ corresponding to the infinitesimal action of $\xi$. 
	
	\item For every smooth function $f$ on $X$, the Hamiltonian vector field $V_f$ is defined by $\iota_{V_f}(\omega)=df$. Explicitly, we have 
	$$
	V_f=\frac{2\pi}{\sqrt{-1}}\left(\frac{\partial f}{\partial z^i}\omega^{i\bar{j}}\frac{\partial}{\partial\bar{z}^j}+\frac{\partial f}{\partial\bar{z}^j}\omega^{\bar{j}i}\frac{\partial}{\partial z^i}\right)
	$$
	(It is easy to check that this vector field indeed satisfies $\iota_{V_f}(\omega)=df$.) 
	And there is
	$$
	\{f,g\}:=-\omega(V_f,V_g).
	$$
	This convention guarantees that the map $f\mapsto V_f$ is a Lie algebra homomorphism:
	$$
	V_{\{f,g\}}=[V_f,V_g].
	$$
	\item On a K\"ahler manifold, under a local complex coordinate system $(z^1,\cdots, z^n)$, we write the symplectic form as 
	$$
	\omega=\frac{\sqrt{-1}}{2\pi}\cdot\omega_{i\bar{j}}dz^i\wedge d\bar{z}^j.
	$$
	We will often use the inverse to the matrix $(\omega_{i\bar{j}})$, which we denote by $(\omega^{i\bar{j}})$. For the entries of this inverse matrix, we have the following conventions:
	$$
	\omega^{i\bar{j}}=-\omega^{\bar{j}i},\hspace{5mm}\omega_{i\bar{j}}\omega^{\bar{j}k}=\delta_i^k
	$$ 
	The Laplacian on $X$ is given by 
	$$
	\Delta=4\pi\cdot\omega^{\bar{j}i}\frac{\partial}{\partial\bar{z}^j}\frac{\partial}{\partial z^i}
	$$
	(For $X=\C^n$, this is the positive Laplacian under the standard metric.)
	\item We also have the following convention for the pre-quantum line bundle on K\"ahler manifolds. The curvature of the pre-quantum line bundle $L$ is 
	$$
	R_L=\omega_{i\bar{j}}dz^i\wedge d\bar{z}^j.
	$$
    For instance, let $e_L$ denote a local holomorphic frame of $L$, with the hermitian inner product 
    $$
    \langle e_L, e_L\rangle=e^{\rho}.
    $$
    Then the function $\rho$ satisfies 
    $$
    \frac{\partial^2\rho}{\partial\bar{z}^j\partial z^i}=-\omega_{i\bar{j}}.
    $$
    So we know that the connection of $L$ is explicitly:
    $$
    \nabla_L(e_L)=\partial\rho\otimes e_L. 
    $$
	\item We use the Einstein summation convention throughout this paper. 
\end{itemize}

\section*{Acknowledgement}
N. C. Leung was supported by grants of the Hong Kong Research Grants Council (Project No. CUHK 14306322 \& CUHK14305923) and direct grants from CUHK. Q. Li was supported by Guangdong Basic and Applied Basic Research Foundation (Project No. 2020A1515011220) and National Science Foundation of China (Project No. 12071204 and Project No. 12471061). Z. N. Ma was supported by National Science Foundation of China (Project No. 1237011981)

\section{Fedosov quantization}\label{section: Fedosov-quantization}

\subsection{Fedosov deformation quantization on K\"ahler manifolds}
\

Recall that a {\em deformation quantization} of a symplectic manifold $(X, \omega)$ is a formal deformation of the commutative algebra $(C^{\infty }(X),\cdot)$ equipped with pointwise multiplication to a noncommutative one $( C^{\infty }( X)[[\hbar]] ,\star) $ equipped with a {\em star product} of the following form
$$
f\star g=fg+\sum_{i\geq 1}\hbar^i\cdot C_i(f,g),
$$
where each $C_i(-,-)$ is a bi-differential operator, so that the leading order of noncommutativity is a constant multiple of the Poisson bracket $\{-,-\}$ associated to $\omega$, i.e.,
\begin{equation}\label{equation: Poisson-bracket}
	C_1(f,g)-C_1(g,f)=\frac{d}{d\hbar}\left( f\star_{\hbar}g - g\star_{\hbar}f\right)\Big|_{\hbar=0}
	=\frac{\sqrt{-1}}{2\pi}\cdot\left\{ f,g\right\}.
\end{equation}

A deformation quantization on a K\"ahler manifold $X$ is called of {\em Wick type} (also known as separation of variables) if all the bi-differential operators $C_l(f,g)$ take holomorphic and anti-holomorphic derivatives of $f$ and $g$ respectively. It is shown in \cite{Karabegov96} that to every Wick type star product, there is an associated closed formal $(1,1)$-form $-\frac{1}{\hbar}\omega+\alpha=-\frac{1}{\hbar} \omega +\alpha_0+ \alpha_1 \hbar + \alpha_2 \hbar^2 + \alpha_3 \hbar^3 + \cdots$ known as the {\em Karabegov form}, which gives rise to a one-to-one correspondence. In this paper, we will focus on Wick type star products on K\"ahler manifolds.
%\begin{example}
%	The Berezin-Toeplitz quantization is a Wick type deformation quantization whose Karabegov form is $-\frac{1}{\hslash} \omega-i\cdot\text{Ric}_X$. 
%\end{example}

In \cite{CLL}, it is shown that every Wick type deformation quantization on a K\"ahler manifold $X$ can be obtained by a Fedosov connection induced from the quantization of $L_\infty$ structure on $X$. 
We now briefly review the construction of this type of Fedosov quantization on K\"ahler manifolds in \cite{CLL}.  %The K\"ahler form on $X$ will always be written in local coordinates as$\omega=\omega_{i\bar{j}}dz^i\wedge d\bar{z}^j$, where we adopt the convention that $\omega^{\bar{k}i}\omega_{i\bar{j}}=\delta_{\bar{j}}^{\bar{k}}$.
We consider the following \emph{Weyl bundles} on $X$:
\begin{align*}\label{equation: Weyl-bundle}
	\W_{X}& := \widehat{\Sym}T^*X, \quad \overline{\W}_X:=\widehat{\Sym}\overline{T^*X}, \quad \W_{X,\C} := \W_{X}\otimes_{\mathcal{C}^\infty_X}\overline{\W}_X=\widehat{\Sym}T^*X_{\C},
\end{align*}
which we call holomorphic, anti-holomorphic and complexified Weyl bundles respectively. 
To give explicit expressions of these bundles, we let $(z^1,\cdots, z^n)$ be a local holomorphic coordinate system on $X$, use $dz^i,d\bar{z}^j$'s to denote $1$-forms in $\A_X^\bullet$ and use $y^i,\bar{y}^j$ to denote sections in $\W_{X,\C}$. The K\"ahler form enables us to define a non-commutative fiberwise Wick product on $\W_{X,\C}$: 
\begin{equation}\label{equation: fiberwise-Wick-product}
	a\star b := \sum_{k\geq 0}\frac{\hbar^k}{k!}\cdot\omega^{i_1\bar{j}_1}\cdots\omega^{i_k\bar{j}_k}\cdot\frac{\partial^k a}{\partial y^{i_1}\cdots\partial y^{i_k}}\frac{\partial^k b}{\partial \bar{y}^{j_1}\cdots\partial \bar{y}^{j_k}}.
\end{equation}

%Throughout this paper, we denote by $\nabla$ the Levi-Civita connection on $X$, and its natural extension to the Weyl bundle $\W_{X,\C}$. By \cite[Proposition 4.1]{Bordemann}, its curvature can be written as a bracket:
%$$\nabla^2=\frac{1}{\hbar}[R_\nabla,-]_\star, $$
%where $R_\nabla=R_{i\bar{j}k\bar{l}}dz^i\wedge d\bar{z}^j\otimes y^k\bar{y}^l\in\A_X^2(\W_{X,\C})$.

 %See \cite[Proposition]{Fed} for this expression of the curvature of $\nabla$; here we are writing it in complex coordinates.

%A natural filtration on these Weyl bundles is defined by polynomial degrees. For instance, $(\overline{\W_X})_{\leq N}$ denotes the sum of anti-holomorphic monomials of polynomial degree $\leq N$. 
The {\em symbol map}
\begin{equation}\label{equation: symbol-map}
	\sigma: \A_X^\bullet(\W_{X,\C})[[\hbar]]\rightarrow\A_X^\bullet[[\hbar]].
\end{equation}
is defined by setting all $y^i,\bar{y}^j$'s to zero. 
%Here $\A_X^\bullet(\W_{X,\C})[[\hbar]]$ denotes the complex of differential forms on $X$ with values in the Weyl bundle. 
On $\A_X^\bullet(\W_{X,\C})[[\hbar]]$, we can define the fiberwise de Rham differential which acts on generators of $\W_{X,\C}$ as $\delta(y^i)=dz^i,\ \delta(\bar{y}^j)=d\bar{z}^j$. We can decompose it as the sum of its $(1,0)$ and $(0,1)$ components $\delta=\delta^{1,0}+\delta^{0,1}$ with respect to the complex structure.

\begin{theorem}[Theorems 2.17 and 2.25 in \cite{CLL}]\label{theorem: Fedosov-connection}
	Let $\alpha = \sum_{i\geq 0}\hbar^i\alpha_i$ be a closed $2$-form on $X$ of type $(1,1)$. Then there exists $I_\alpha = I+ J_\alpha \in \mathcal{A}_X^{0,1}(\mathcal{W}_{X,\mathbb{C}})[[\hbar]]$, such that the connection 
	%\begin{equation}\label{equation: Fedosov-equation}
	%	\nabla I_\alpha - \delta I_\alpha + \frac{1}{\hbar} I_\alpha\star I_\alpha + R_\nabla=\alpha.
	%\end{equation}
	\begin{equation}\label{equation: connection-D-alpha}
	D_{\alpha} := \nabla-\delta+\frac{1}{\hbar}[I_\alpha, -]_{\star}
	\end{equation}
	is a Fedosov abelian connection. 
	The deformation quantization associated to the flat connection $D_{\alpha}$ is a Wick type star product with Karabegov form $K_{\alpha}$ given by $-\frac{1}{\hbar}\cdot\omega+\alpha$. 
\end{theorem}
   For an explicit expression of the term $I_\alpha$, we refer to \cite{CLL}. Since the operator $\delta$ can also be written as a bracket (with respect to the fiberwise Wick product $\star$):
   $$
   \delta=\frac{1}{\hbar}\omega_{i\bar{j}}[dz^i\otimes\bar{y}^j-d\bar{z}^j\otimes y^i,-]_\star.
   $$
   the connection $D_\alpha$ can also be written as 
	$D_{\alpha}=\nabla+\frac{1}{\hbar}[\gamma_\alpha,-]_\star$.
	The flatness of $D_{\alpha}$ is then equivalent to the following {\em Fedosov equation}:
	\begin{equation}\label{equation: Fedosov-equation-gamma}
		\nabla\gamma_\alpha+\frac{1}{\hbar}\gamma_\alpha\star\gamma_\alpha+R_\nabla=\frac{2\pi}{\sqrt{-1}}\left(-\frac{1}{\hbar}\cdot\omega+\alpha\right).
	\end{equation}
	
The main result in \cite{Fed} is that there is a one-to-one correspondence between formal smooth functions on $X$ and flat sections of the Weyl bundle under Fedosov's abelian connection via the symbol map, which induces a star product on formal smooth functions.

\subsection{Quantizable functions}
\

%The notion of {\em canonical quantization} is roughly speaking the procedure of quantizing functions to operators on Hilbert spaces.  
In general, we only need to quantize a sub-class of smooth functions known as {\em quantizable functions}.  In \cite{CLL23}, quantizable functions on pre-quantizable K\"ahler manifolds are defined by imposing certain finiteness conditions on the associated flat section in the Weyl bundle:
%We only need a sub-class of smooth functions which can be quantized, which are called the {\em quantizable functions}.  This sub-class should satisfy several properties: it should be closed under the quantum product and there should be enough such functions. A typical example of quantizable functions is the polynomial functions on $\C^n$. 
%To define quantizable functions on general K\"ahler manifolds, we need the following
\begin{lem-defn}\label{definition:degree_1_quantizable_function}
	%We define a \emph{weight} on $\W_{X,\C}[[\hbar]]$ by assigning weights on its generators:
	%\begin{equation}\label{equation: weights-formal-Weyl-bundle}
	%	|y^i|=0, \hspace{2mm} |\bar{y}^j|=2, \hspace{2mm}|\hbar|=2.
	%\end{equation}
	%The fiberwise Wick product $\star$ preserves this weight. 
	Suppose $\alpha$ is a polynomial in $\hbar$, then the subspace  $$\A_X^\bullet\left(\Sym^\bullet\overline{T^*X}\otimes\W_X[\hbar]\right)\subset\A_X^\bullet(\W_{X,\C}[[\hbar]]),$$
	consisting of sections which are polynomials in $\bar{y}_i$'s and $\hbar$, is closed under the connection $D_\alpha$. 
	A formal function $f\in C^\infty(X)[[\hbar]]$ is called {\em quantizable} if its associated flat section $O_f$ lives in $\Sym^\bullet\overline{T^*X}\otimes\W_X[\hbar]$. These formal functions are closed under the associated star product defined via $D_\alpha$. 
\end{lem-defn}

From now on, we assume that $\alpha$ is a polynomial in $\hbar$. We can take the evaluation $\hbar=1/k$ in $D_\alpha$ to obtain a flat connection $D_{\alpha,k}$ (in the formula of the fiberwise formal product $\star$, $\hbar$ should also be evaluated $\hbar=1/k$), which preserves the space $\left(\A_X^\bullet(\Sym^\bullet\overline{T^*X}\otimes\W_X), D_{\alpha,k}\right)$. Similarly, we can define non-formal {\em quantizable functions} of level $k$ whose associated flat section $O_f$ lives in this subspace.  
\begin{notation}
	There are different types of star products, and we list here the notations for readers convenience. We let $\star$ denote the formal fiberwise Wick product on $\W_{X,\C}$, and let $\star_k$ denote the (non-formal) Wick product with evaluation $\hbar=1/k$. Given a Karabegov form $K=-\frac{1}{\hbar}\omega+\alpha$, we let $\star_\alpha$ denote the corresponding Wick type deformation quantization on $C^\infty(X)[[\hbar]]$. In particular, we let $\star_{BT}$ denote the Berezin-Toeplitz quantization on $X$.
\end{notation}

%Suppose $\alpha$ is a polynomial in $\hbar$, the connection $D_\alpha$ (and also its evaluations $D_{\alpha,k}$) increase the polynomial degrees in $\Sym^\bullet\overline{T^*X}$ by a fixed number, we obtain the following sub-complex:
%$$\left(\A_X^\bullet(\Sym^\bullet\overline{T^*X}\otimes\W_X), D_{\alpha,k}\right),$$
%which is closed under the fiberwise star product $\star_k$. 
\begin{definition}\label{definition: quantizable functions}
	A flat section of $\Sym^\bullet\overline{T^*X}\otimes \W_X$ under the Fedosov connection $D_{\alpha,k}$ is called a {\em (non-formal) quantizable function of level $k$}. These sections form a sheaf of algebras on $X$ under the product $\star_k$ on $\W_{X,\C}$.  We denote this sheaf and its global sections by $\mathcal{C}_{\alpha,k}^\infty$ and $C^\infty_{\alpha,k}(X)$ respectively.  
\end{definition}

Given a flat section $\gamma$ of $D_{\alpha,k}$, we can take its symbol $\sigma(\gamma)$ to get an element in $C^{\infty}(X)$. By \cite[Proposition 2.26.]{CLL23}, we learn that $\gamma$ is uniquely determined by $\sigma(\gamma)$ for $k>>0$ and we can call $\sigma(\gamma)$ a quantizable function of level $k$ by abuse of notation. Given a formal quantizable function $f$ with the associated flat section $O_f$, we obtain a non-formal quantizable function of level $k$ by taking the evaluation $O_f$ at $\hbar=1/k$. 
\begin{definition}\label{definition:degree_in_weyl}
	On the formal Weyl bundle $\W_{X,\C}[[h]]$, we can define a weight by assigning weights on generators: $|y^i|=0,  |\bar{y}^j|=|\hbar|=1$. We let $(\W_{X,\C}[[\hbar]])_N$ denote sums of monomials with weights $\leq N$. Similarly we can define a weight on $\mathcal{C}_{\alpha,k}^\infty$ by assigning $|y^i|=0$ and $|\bar{y}^j| = 1$ counting polynomial degree in $\bar{y}_j$'s. And we can correspondingly define {\em degree $N$} (formal or non-formal) quantizable functions. 
\end{definition}

It is easy to see that for every $k\in\C\setminus\{0\}$, the sheaf $\mathcal{C}_{\alpha,k}^\infty$ is closed under the star product $\star_k$ defined via the Fedosov connection $D_{\alpha,k}$. This star product is compatible with the above filtration in the sense that $(\mathcal{C}_{\alpha,k}^\infty)_{N_1}\star_k (\mathcal{C}_{\alpha,k}^\infty)_{N_2}\subset (\mathcal{C}_{\alpha,k}^\infty)_{N_1+N_2}$. 

\subsection{Hilbert spaces and differential operators via Fedosov quantization}\label{section: Hilbert-spaces-quantizable-functions}
\

The K\"ahler form on $X$ enables us to define the fiberwise Bargmann-Fock action: a monomial in $\W_{X,\C}$ acts as a differential operator on $\W_X$ as
\begin{equation}\label{equation: Wick-ordering-formula}
	y^{i_1}\cdots y^{i_k}\bar{y}^{j_1}\cdots\bar{y}^{j_l}\mapsto \left(-\hbar\right)^l\omega^{p_1\bar{j}_1}\cdots\omega^{p_l\bar{j}_l}\frac{\partial}{\partial y^{p_1}}\circ\cdots  \frac{\partial}{\partial y^{p_l}}\circ m_{y^{i_1}\cdots y^{i_k}}.
\end{equation}
This action is compatible with the fiberwise (formal) Wick product on the Weyl bundle $\mathcal{W}_{X,\C}$ in the following sense: let $\alpha_1,\alpha_2$ be sections of $\W_{X,\C}$, and let $s$ be a section of $\W_X$:
$$
(\alpha_1\star_\hbar\alpha_2)\circledast s=\alpha_1\circledast(\alpha_2\circledast s),
$$
making the holomorphic Weyl bundle $\W_X$ a sheaf of $\W_{X,\C}$-modules. 
In the rest of this section, we will mainly focus on the Berezin-Toeplitz deformation quantization, which is of Wick type with $\alpha_{\text{BT}} =  Ric_X$ or the Karabegov form as
\begin{equation}\label{equation: Karabegov-form-Berezin-Toeplitz}
K_{\text{BT}} = -\frac{1}{\hbar}\cdot\omega+Ric_X.
\end{equation}
%We will drop the subscript $\text{BT}$ in our notation for simplicity. 
\begin{definition}
	For every positive integer $k$,  we define the \emph{level $k$ Bargmann-Fock sheaf} $\mathcal{F}_{L^{\otimes k}}$ by twisting $\W_X$ with the $k$-th tensor power of the prequantum line bundle $L$:
	$$
	\mathcal{F}_{L^{\otimes k}}:=\mathcal{W}_{X}\otimes_{\mathcal{O}_X}L^{\otimes k}.
	$$
\end{definition}
The Bargmann-Fock sheaf forms a module over the bundle of finite order Weyl bundle via the fiberwise Bargmann-Fock action. It is shown in \cite{CLL} that there also exists a Fedosov type abelian connection $D_{BT,k}$ on $\mathcal{F}_{L^{\otimes k}}$ given by
$$
D_{BT,k} = \nabla + k\cdot  \gamma_{BT}  \circledast_k. 
$$ Similar to the case in the algebra of quantizable function of level $k$, we also have the symbol map giving an isomorphism from flat sections of the Bargmann-Fock sheaf to holomorphic sections of $L^{\otimes k}$:
\begin{equation}\label{equation: symbol-map-Bargmann-Fock}
	\sigma: \mathcal{F}_{L^{\otimes k}}\rightarrow L^{\otimes k}.
\end{equation}
\begin{remark}
	Throughout this paper, when we say flat sections, we always refer to flat sections of either the Weyl or the Bargmann-Fock bundle with respect to the corresponding Fedosov connection. 
\end{remark}
For later computation, we will need an explicit expression of a flat section of $\mathcal{F}_{L^{\otimes k}}$ associated to a holomorphic section $s\in H^0(X, L^{\otimes k})$. Let $e_{L^{k}}$ denote a local holomorphic frame of $L^{\otimes k}$. Suppose the hermitian metric of $L^{\otimes k}$ is given locally by $\langle e_{L^{k}},e_{L^{k}}\rangle=e^{k\cdot\rho}$. 
%The connection $\nabla_{L^{\otimes k}}$ can then be written explicitly as
%$$\nabla_{L^{\otimes k}}(e_{L^{k}})=k\cdot\partial\rho\otimes e_{L^{k}},$$
%and the prequantum condition implies that $\bar{\partial}\partial\rho=\omega$.  
We define a local section of  $\mathcal{F}_{L^{\otimes k}}$ by $\beta:=\sum_{l\geq 1}(\tilde{\nabla}^{1,0})^l(\rho)$, which is shown in \cite{CLL23} to be flat under $D_{BT,k}$ such that $\sigma(\beta)=e_{L^k}$. For general $s\in H^0_{\bar{\partial}}(X,L^{\otimes k})$, it must be locally of the form $g\cdot e_{L^k}$ for a holomorphic function $g$.  Then the flat section corresponding to $s$ is locally given by 
\begin{equation}\label{equation: flat-section-Hilbert-spaces}
	O_g\cdot e^{k\cdot\beta}\otimes e_{L^{k}}.
\end{equation}

Quantizable functions on the Bargmann-Fock sheaves via the formula in equation \eqref{equation: Wick-ordering-formula}. 
Since the Fedosov connections $D_{BT,k}$ on the Weyl bundle and the Bargmann-Fock sheaves are compatible in the sense that
\begin{equation}\label{equation: Fedosov-connection-compatibility}
	D_{BT,k}(\alpha_1\circledast_k s)=D_{BT,k}(\alpha_1)\circledast_k s+\alpha_1\circledast_k D_{BT,k}(s).
\end{equation}
(This is also the reason we abuse the notation $D_{BT,k}$). It follows that the output of the action of a quantizable function on a flat section of $\mathcal{F}_{L^{\otimes k}}$ is still flat. This action is obviously local, thus quantizable functions act as holomorphic differential operators. And we obtain a homomoprhism of sheaf of algebras:
\begin{equation}\label{equation: isomorphism-quantizable-function-differential-operator}
\varphi: \mathcal{C}_{BT,k}^\infty\rightarrow \mathcal{D}(L^{\otimes k}).
\end{equation}
It is shown in \cite{CLL23} that the homomorphism $\varphi$ is actually an isomorphism. 
For every quantizable function $f$ with highest weight term given by $f_{j_1\cdots j_m}\bar{y}^{j_1}\cdots \bar{y}^{j_m}$, the principal symbol $\varphi(f)$ of the corresponding differential operator is given by
$$
\varphi(f) = \left(-\frac{1}{k}\right)^m\cdot f_{j_1\cdots j_m}\omega^{i_1\bar{j}_1}\cdots\omega^{i_m\bar{j}_m}\frac{\partial}{\partial z^{i_1}}\cdots\frac{\partial}{\partial z^{i_m}}.
$$

%\begin{theorem}\label{theorem: Bargmann-Fock-isomorphic-prequantum}
%	Suppose that $X$ is a K\"ahler manifold equipped with a prequantum line bundle $L$, and we choose the formal closed $(1,1)$-form $\alpha$ as in equation \eqref{equation: Karabegov-form-Berezin-Toeplitz}. Then for any positive integer $k$, the symbol map gives a sheaf isomorphism from $\mathcal{F}_{L^{\otimes k}}$ under the connection $D_{\alpha,k}$ to the sheaf of holomorphic sections of $L^{\otimes k}$.
%\end{theorem}

\section{Degree $1$ quantizable functions in Berezin-Toeplitz quantization}\label{section: first-order-quantizable-function}

%We have seen in the previous section that the quantum Hamiltonian associated to an infinitesimal holomorphic isometry of a K\"ahler manifold gives rise to a degree $1$ quantizable function by taking the evaluation $\hbar=1/k$. 

In this section, we will give a full description of (both formal and non-formal) quantizable functions in Berezin-Toeplitz deformation quantization of at most degree $1$ in the anti-holomorphic Weyl bundle $\Sym^\bullet\overline{T^*X}$. We will simply call these degree $1$ quantizable function. We begin with a full description of formal degree $1$ quantizable functions.

%There are two main results in this section: first we describe all degree $1$ formal quantizable functions, next we will show that this also implies the result on non-formal degree $1$ quantizable functions, since every such function can be obtained by taking the evaluation of a formal quantizable function. 

\subsection{Degree $1$ formal quantizable functions}\label{section: formal-degree-1-functions}
\

%In this section, we will give a full description of degree $1$ formal quantizable functions in Berezin-Toeplitz quantization. 
\begin{prop}\label{proposition: degree-1-formal-quantizable-function}
	Suppose $f_0\in C^\infty(X)$ is a smooth function such that 
	\begin{equation}\label{equation: f-0}
	\nabla^{0,1}\left(\frac{\partial f_0}{\partial\bar{z}^j}\bar{y}^j\right)=0.
	\end{equation}
	Then $f_0-\hbar(\frac{1}{4\pi}\Delta f_0)$ is a formal degree $1$ quantizable function. Conversely, suppose $f_0-\hbar f_1$ is a formal degree $1$ quantizable function, then equation \eqref{equation: f-0} is satisfied and $f_1-\frac{1}{4\pi}\Delta f_0$ is a holomorphic function. Thus if $X$ is compact, then every formal degree $1$ quantizable function can be written as $f=f_0-\hbar\cdot(\frac{1}{4\pi}\Delta f_0+c)$ for some constant $c\in\C[\hbar]$, where $f_0$ satisfies equation \eqref{equation: f-0}.
\end{prop}
\begin{proof}
	From \cite[Proposition 2.19.]{CLL} we see that $O_f$ is determined by the iterative equation:
	\begin{equation}\label{equation:Of_iteration}
	O_f = f + \delta^{-1}(\nabla O_f + \frac{1}{\hbar}[I_\alpha,O_f]_{\star}).
	\end{equation}
	By making use of \cite[Proposition 2.24.]{CLL3}, which also holds for formal quantizable function, we can simply restrict ourselves to components lying in $\overline{\mathcal{W}}_X[\hbar]$. We consider the components $\sum_{i \geq 0}\beta_{j,i} \hbar^{i}$ of the iteration lying in $\Sym^j(\overline{T^*X})[\hbar]$ according to filtration by weight, the term of weight $1$ is $\beta_{1,0} = \frac{\partial f_0}{\partial\bar{z}^j}\bar{y}^j$ and $\nabla^{0,1}\beta_{1,0}   = 0$ from which follows from equation \eqref{equation: f-0}, and is equivalent to that 
	$$
	\frac{\partial f_0}{\partial\bar{z}^j}\omega^{i\bar{j}}
	$$
	being holomorphic for all indices $i$. We claim that the next term in the iteration vanishes, i.e. 
	$$
\left[ (Ric_X)_{k\bar{l}}d\bar{z}^l\otimes y^k, \frac{\partial f_0}{\partial\bar{z}^j}\bar{y}^j\right]_{\star}  + \frac{\hbar}{4\pi} \nabla^{0,1}(\Delta f_0) = 0,
	$$ 
	and hence we see that $\beta_{j,i} = 0$ for $j>1$ or $i>0$. 
	First, we have the following computation
	\begin{align*}
		\left[(Ric_X)_{k\bar{l}}d\bar{z}^l\otimes y^k, \frac{\partial f_0}{\partial\bar{z}^j}\bar{y}^j\right]_{\star}=&\hbar\cdot\frac{\partial f_0}{\partial\bar{z}^j}\cdot\omega^{k\bar{j}}\cdot (Ric_X)_{k\bar{l}}d\bar{z}^l\\
		=&-\hbar\cdot\frac{\partial f_0}{\partial\bar{z}^j}\cdot\omega^{k\bar{j}}\cdot \frac{\partial\Gamma^i_{ik}}{\partial\bar{z}^l}d\bar{z}^l. 
	\end{align*}
	From the equation \eqref{equation: f-0} and taking a local normal coordinate \cite[Definition 2.8.]{CLL3} centered at an arbitrary point $x_0\in X$, under which there is $\omega_{i\bar{j}}(z)=\delta_{ij}+O(|z|^2)$. And we have
	\begin{align*}
		0=\frac{\partial^2}{\partial\bar{z}^l\partial z^i}\left(\frac{\partial f_0}{\partial\bar{z}^j}\omega^{i\bar{j}}\right)\bigg|_{x_0}=&\left(\frac{\partial^2}{\partial\bar{z}^l\partial z^i}\left(\frac{\partial f_0}{\partial\bar{z}^j}\right)\omega^{i\bar{j}}+\frac{\partial f_0}{\partial\bar{z}^j}\cdot\frac{\partial^2\omega^{i\bar{j}}}{\partial\bar{z}^l\partial z^i}\right)\bigg|_{x_0}\\
		=&\left(-\frac{1}{4\pi}\frac{\partial\Delta f_0}{\partial\bar{z}^l}+\frac{\partial f_0}{\partial\bar{z}^j}\cdot\frac{\partial^2\omega^{i\bar{j}}}{\partial\bar{z}^l\partial z^i}\right)\bigg|_{x_0} \\
		=& \left(-\frac{1}{4\pi}\frac{\partial\Delta f_0}{\partial\bar{z}^l}+\frac{\partial f_0}{\partial\bar{z}^j}\cdot 	\omega^{k\bar{j}}\cdot \frac{\partial\Gamma^i_{ik}}{\partial\bar{z}^l}\right)\bigg|_{x_0}
	\end{align*}
	where the last equality follows from
	$$
	\left(\omega^{k\bar{j}}\cdot \frac{\partial\Gamma^i_{ik}}{\partial\bar{z}^l}\right)\bigg|_{x_0}=\left(\frac{\partial^2\omega^{i\bar{j}}}{\partial\bar{z}^l\partial z^i}\right)\bigg|_{x_0}.
	$$
	Putting these equations together we obtain the desired equation.  
	
	Conversely, if $f_0 - \hbar f_1$ is a formal degree $1$ quantizable function we see that $\beta_{1,0} =  \frac{\partial f_0}{\partial\bar{z}^j}\bar{y}^j$ which must satisfy $\nabla^{0,1}(\beta_{1,0}) = 0$. Therefore $\hbar (f_1 - \frac{1}{4\pi}\Delta f_0) $ is again a formal degree $1$ quantizable function, which forces $\frac{\partial }{\partial\bar{z}^j} (f_1 - \frac{1}{4\pi}\Delta f_0) = 0$. 
	%We would like to compare it with the following:
	%\begin{align*}
	%	\frac{\partial \Delta f_0}{\partial\bar{z}^l}=&\frac{\partial^3 f_0}{\partial\bar{z}^l\partial\bar{z}^j\partial z^i}\omega^{i\bar{j}}
	%\end{align*}
\end{proof}

\subsection{Surjectivity of the evaluation map}
\

We have seen that given a formal quantizable function $f$, we can take the evaluation of $O_f$ at $\hbar=1/k$ for any positive integer $k$ to obtain a non-formal quantizable function. A natural question is if the converse is true. It is shown in \cite{CLL23} locally this is true, since sheaf-theoretically, we can identify quantizable functions with holomorphic differential operators. We will show that this is also true for global degree $1$ quantizable functions. Equivalently, every degree $1$ non-formal quantizable function is essentially formal:

\begin{theorem}\label{theorem: every-degree-1-quantizable-function-is-formal}
For the Berezin-Toeplitz deformation quantization, every non-formal degree $1$ quantizable function is the evaluation of a formal degree $1$ quantizable function. 
\end{theorem}
\begin{proof}
Suppose $\beta$ is the flat section corresponding to a degree one quantizable function of level $k$. Let $\beta_{0,1}$ denote the component in $\beta$ of polynomial degree $0$ and $1$ in $\W_X$ and $\overline{\W}_X$ respectively, which we can write as 
$$
\beta_{0,1}=(\beta_{0,1})_{\bar{j}}\bar{y}^j. 
$$ 
It follows from the flatness $D_{BT,k}(\beta)=0$ that 
$$
\nabla^{0,1}(\beta_{0,1})=\nabla^{0,1}((\beta_{0,1})_{\bar{j}}\bar{y}^j)=0. 
$$
We must also have 
$$
\bar{\partial}\left((\beta_{0,1})_{\bar{j}}d\bar{z}^j\right)=0.
$$
It follows from the $\bar{\partial}$-Poincare lemma that on a local open set $U_\gamma\subset X$, there exists a function $f_\gamma$, such that $\bar{\partial}f_\gamma=(\beta_{0,1})_{\bar{j}}d\bar{z}^j$. By Proposition \ref{proposition: degree-1-formal-quantizable-function}, on each $U_\gamma$, 
$$
f_\gamma-\frac{\hbar}{4\pi}\Delta f_\gamma
$$
is a formal degree $1$ quantizable function.
 Let $\eta_{\gamma,k}=(O_{f_\gamma-1/{4\pi\hbar\Delta f_\gamma}})|_{\hbar=1/k}$ be the quantizable function on $U_\gamma$ corresponding to $f_\gamma-\frac{1}{4\pi k}\Delta f_\gamma$ evaluated at $\hbar=1/k$. Then $\eta_{\gamma,k}$ and $\beta|_{U_\gamma}$ have the same highest degree term $(\beta_{0,1})_{\bar{j}}\bar{y}^j$, and their difference must be the flat section $O_{h_\gamma}$ corresponding to a holomorphic function $h_\gamma\in\mathcal{O}_X(U_\gamma)$. On the overlap $U_\gamma\cap U_\xi$, there is 
$$
(\eta_{\gamma,k}-\eta_{\xi,k})|_{U_\gamma\cap U_\xi}=(\eta_{\gamma,k}-\beta)|_{U_\gamma\cap U_\xi}-(\eta_{\xi,k}-\beta)|_{U_\gamma\cap U_\xi}=(O_{h_\gamma})|_{U_\gamma\cap U_\xi}-(O_{h_\xi})|_{U_\gamma\cap U_\xi}.
$$
It is then easy to see that the functions 
$$
f_\gamma-h_\gamma-\frac{\hbar}{4\pi}\Delta f_\gamma
$$
glue to a global degree $1$ formal quantizable function, whose evaluation at $\hbar=1/k$ is exactly $\beta$. 
\end{proof}

\section{Quantum Hamiltonian functions as degree $1$ quantizable functions}\label{section: quantum-Hamiltonian}
In this section, we will give a class of examples of degree $1$ formal quantizable functions which arises from symmetries of deformation quantizations. We call these functions {\em quantum Hamiltonians}. We will also show that a degree $1$ formal quantizable function is a quantum Hamiltonian if and only if it is real. 
%Recall that a {\em quantum moment map} is a Lie algebra homomorphism $\mu_\hbar^*:\mathfrak{g}\rightarrow \left(C^\infty(X)[[\hbar]],\frac{i}{\hslash}[\cdot,\cdot]_\star\right)$ , such that for any $v\in\mathfrak{g}$:
%$$v(f)=\frac{i}{\hslash}[\mu_\hbar^*(v),f]_\star. $$

%In this subsection, we will describe the flat sections associated to the quantum moment maps.

%We can define a similar and closely related notion of quantum Hamiltonian:
\subsection{Quantum Hamiltonian}
\

For this section, we work with a general Wick type deformation quantization whose corresponding Karabegov form is $-\frac{1}{\hbar}\omega+\alpha$. Let us recall that a vector field $V$ is called Hamiltonian if there exists a smooth function $f$ such that,  for any smooth function $g\in C^\infty(X)$, there is
$$
V(g)=\{f,g\}.
$$
Then $f$ is called a (classical) Hamiltonian function of the vector field $V$. In particular, the vector field $V$ must preserve $\omega$, i..e, it is symplectic. 

We have a quantum version of notion of Hamiltonian. Here we refer to \cite{Muller-Neumaier, Muller-Neumaier-2} for the notion and results. We first start from the following definition:
\begin{definition}
	Let $V$ be a vector field on $X$. We say that $V$ is quasi-inner, if there exists a smooth formal function $f\in C^\infty(X)[[\hbar]]$, such that 
	$$
	V(g)=\frac{1}{\hbar}[f,g]_\star. 
	$$
	This formal function $f$ is said to be a quantum Hamiltonian associated to $V$, with respect to this star product $\star$. 
\end{definition}
It is easy to see that if a quantum Hamiltonian exists, it is unique up to formal constant functions. It is also clear from the definition that such a vector field must act as a derivation on formal smooth functions. And we can apply the following criterion on such vector fields (\cite{Muller-Neumaier-2}[Proposition 3.2])
\begin{prop}\label{proposition: vector-field-derivation}
Let $\star$ be a Wick type star product on $X$ with Karabegov form $K$. A vector field $V$ on $X$ is a derivation with respect to $\star$ if and only if 
$$
\mathcal{L}_V(J)=0, \hspace{3mm} \textit{and}\hspace{3mm}\mathcal{L}_V(K)=0.
$$
Here $J$ denotes the complex structure on $X$. 
\end{prop}
Since the first term of the Karabegov form $K$ is the symplectic form $\omega$, it is clear that such a vector field $V$ must be symplectic. Moreover, since $V$ is also compatible with the complex structure, it must be Killing. And we immediately have the following for the Berezin-Toeplitz star product $\star_{\text{BT}}$. 
\begin{corollary}\label{corollary: vector-field-compatible-Berezin-Toeplitz}
	A vector field $V$ on $X$ is compatible with the Berezin-Toeplitz deformation quantization if and only if it is both symplectic and Killing. 
\end{corollary}
\begin{proof}
	The Karabegov form of the Berezin-Toeplitz deformation quantization is 
	$$
	K_{\text{BT}}=-\frac{1}{\hbar}\cdot\omega+ Ric_X.
	$$
	If a symplectic vector field is Killing, then it preserves both the Ricci form  and the symplectic form, and thus the Karabegov form. 
\end{proof}

It is shown in \cite{Muller-Neumaier-2} that for every vector field compatible with a star product, quantum Hamiltonians always exist locally. And there is a cohomological obstruction for their global existence:
\begin{theorem}\label{theorem: symmetry-Kahler}
	Let $\star_\alpha$ be a star product of Wick type on $X$, with Karabegov form $K=-\frac{1}{\hbar}\omega+\alpha$. Suppose $V$ is a vector field on $X$ which acts on $C^\infty(X)[[\hbar]]$ as a derivation with respect to $\star$. Then $V$ is quasi-inner if there exists a formal smooth function $f$ such that 
	\begin{equation}\label{equation: obstruction-quantum-hamiltonian}
	df=-\frac{2\pi\hbar}{\sqrt{-1}}\cdot\iota_VK
	\end{equation}
\end{theorem}

Since different choices of quantum Hamiltonians differ by a constant formal function (valued in $\C[[\hbar]]$), there should be a canonical section of the Weyl bundle modulo the constant terms, in terms of the Fedosov quantization scheme.  We will give an explicit expression of this Weyl bundle sections. In particular, it will be almost obvious from these expressions that when the Karabegov form is {\em admissible}, quantum Hamiltonians are always quantizable functions. 

Recall that the Fedosov connection associated to $\alpha$ is the form 
$$
D_\alpha=\nabla+\frac{1}{\hbar}[\gamma_\alpha,-]_\star. 
$$
Suppose the symplectic vector field $V$ preserves the complex structure and Karabegov form $-\frac{1}{\hbar}\omega+\alpha$. Since all the components in $D_\alpha$ (and also the bracket $[-,-]_\star$) only depend on the complex, symplectic and Riemannian structure, and that the vector field $V$ preserves all these structures, there is the following commutativity relation:
\begin{equation}\label{equation: Lie-derivative-commutes-Fedosov}
	[\mathcal{L}_V, D_\alpha]=0. 
\end{equation}

 It is shown in \cite{CLL23} that there exists a section $\eta=\eta_{i\bar{j}}y^i\bar{y}^j\in\Gamma(X,\W_{X,\C})$ (of polynomial degree $2$ and type $(1,1)$), such that 
\begin{equation}\label{equation: moment-as-bracket}
	\mathcal{L}_{V}-\nabla_{V}=\frac{1}{\hbar}[\eta,-]_\star.
\end{equation}
%We claim that modulo constant terms (i.e., those terms of polynomial degree $0$) in $\W_{X,\C}$, there is a flat section $\beta_V$ is of the form:

%In particular, we have the following terms in $\beta_V$ which is of polynomial degree $1$ in the anti-holomorphic Weyl bundle:
%$$
%\sum_{m\geq 0}(\tilde{\nabla}^{1,0})^m\left(\omega_{i\bar{j}}\cdot \iota_{V} dz^i\otimes \bar{y}^j\right).
%$$
%The proof of the the following lemma is similar to that of \cite[Proposition 5.4]{CLL23}. For the completeness of this paper, we give a proof here:

\begin{lemma}\label{lemma:lie_action_on_flat_section}
	Suppose $V$ is a vector field on $X$ that preserves the complex structure and Karabegov form $K=-\frac{1}{\hbar}\omega+\alpha$. We define a section of the Weyl bundle by
	\begin{equation}\label{equation: flat-section-moment-map}
		\beta_V:=\eta-\iota_V \gamma_\alpha.
	\end{equation}
		For any smooth function $f\in C^\infty(X)$,	 the Lie derivative $\mathcal{L}_V$ action on $O_f$ can be described by the bracket with $\beta_V$:
	$$
	\mathcal{L}_{V}(O_f) = O_{\mathcal{L}_V(f)}=\frac{1}{\hbar}[\beta_V, O_f]_\star
	$$
\end{lemma}
\begin{proof}
	From the construction of $\beta_V$ and equation \eqref{equation: moment-as-bracket}, we have 
	$$
	\frac{1}{\hbar}[\beta_V,-]_\star=\mathcal{L}_V-[D_\alpha,\iota_V].
	$$
	And we have 
	\begin{align*}
		\frac{1}{\hbar}[\beta_V,O_f]_\star=&\mathcal{L}_V(O_f)-[D_\alpha,\iota_V](O_f)=\mathcal{L}_V(O_f)=O_{\mathcal{L}_V(f)}.
	\end{align*}
	In the second equality we have used the flatness of $O_f$. For the third equality, notice that $\mathcal{L}_V(O_f)$ is also flat due to equation \eqref{equation: Lie-derivative-commutes-Fedosov}, and that its symbol is exactly $\mathcal{L}_V(f)$. 
\end{proof}

%\begin{remark}
%	It is clear that  for any $g\in\mathfrak{g}$, the image of $g$ under the quantum moment map is a quantum Hamiltonian of the associated vector field $V_g$. 
%\end{remark}
To summarize, we have the following:
\begin{theorem}
	Let $V$ be a vector field on $X$ as in the previous lemma.  Then the associated Lie derivative is compatible with the Fedosov connection $D_\alpha$:
	\begin{equation*}
		[\mathcal{L}_V, D_\alpha]=0. 
			\end{equation*}
		The action of $V$ on smooth formal functions has a Fedosov lifting in the following sense:
		$$
		\mathcal{L}_V(O_f)=O_{V(f)}.
		$$
		
		 A quantum Hamiltonian $\mu_V$ associated to $V$ always exists locally, whose flat section modulo the constant term is explicitly given by 
		\begin{equation}\label{equation: flat-section-modulo-constant-term}
		\beta_V:=\eta-\iota_V \gamma_\alpha.
		\end{equation}

		The symbol (i.e., the local quantum Hamiltonian) is determined by the following equation which always has solutions locally on $X$:
			\begin{equation}\label{equation: condition-quantum-moment-map}
			-\frac{2\pi\hbar}{\sqrt{-1}}\cdot\iota_VK=d\mu_V.
		\end{equation}
		Global quantum Hamiltonian $\mu_V$ exists globally if the $1$-form $-\frac{2\pi\hbar}{\sqrt{-1}}\cdot\iota_VK$ is exact.
\end{theorem}
\begin{proof}
	%For the first statement, suppose we are given two smooth functions $f,g\in C^\infty(X)$ with $O_f$ and $O_g$  the corresponding flat sections respectively, there is
	%$$D_\alpha\left(\mathcal{L}_V(O_f\star O_g)\right)=0.$$
	%On the other hand, we have
	%\begin{align*}
	%	\mathcal{L}_V(O_f\star O_g)=&\mathcal{L}_V(O_f)\star O_g+O_f\star\mathcal{L}_V(O_g).
	%\end{align*}
	%By looking at the symbol of two sides of the above equation, we obtain equation \eqref{equation: vector-field-derivation}.
	The first two statements follow from the previous discussions. We only prove the remaining statements here: using equations \eqref{equation: moment-as-bracket} and \eqref{equation: flat-section-moment-map}, we have
	\begin{align*}
		\frac{1}{\hbar}[[\beta_V,-]_\star, D_\alpha]=&\left([\mathcal{L}_V,D_\alpha]-[[D_\alpha,\iota_V], D_\alpha]\right)=-\left([[D_\alpha,\iota_V], D_\alpha]\right)\\
		=&-\frac{1}{2}\iota_V\left([D_\alpha, D_\alpha]\right)=0;
	\end{align*}
	here we used equation \eqref{equation: Lie-derivative-commutes-Fedosov} in the above lemma in the second equality. On the other hand, there is
	$$
	\frac{1}{\hbar}[[\beta_V,-]_\star, D_\alpha]=\pm\frac{1}{\hbar}[D_\alpha(\beta_V),-]_\star=\pm\frac{1}{\hbar}[D_\alpha(\beta_V),-]_\star=0. 
	$$
	It follows that  $D_\alpha(\beta_V)\in\A^1(X)$  (since $\A^\bullet(X)$ is the center of $\A_X^\bullet(\W_{X,\C})$). Moreover, these $1$-forms must be closed since $D_\alpha^2(\beta_V)=d_X(D_\alpha(\beta_V))=0$. Recall that we have $\gamma_{\alpha} = I+J_{\alpha} - \omega_{i \bar{j}}(dz^i \otimes \bar{y}^{j} - d\bar{z}^{j} \otimes y^i)$. We have the terms 
	$$
	-\left(-\hbar\cdot\alpha_{i\bar{j}} \iota_V (d\bar{z}^j) \otimes y^i - \omega_{i \bar{j}} \iota_{V}(dz^i \otimes \bar{y}^{j} - d\bar{z}^{j} \otimes y^i)\right)
	$$
	in the expression of $\beta_V$ which contains only polynomial of degree $1$ in $y^i$'s or $\bar{y}^j$'s. It turns out that the components of $D_{\alpha}(\beta_V)$ lying in $\A^1(X)$ are contributed from the following terms
	\begin{align*}
		&D_{\alpha}(\beta_V) \\
		=&  - \delta	\left(\hbar\cdot\alpha_{i\bar{j}} \iota_V (d\bar{z}^j) \otimes y^i  + \omega_{i \bar{j}} \iota_{V}(dz^i \otimes \bar{y}^{j} - d\bar{z}^{j} \otimes y^i)\right) + \frac{1}{\hbar}  \left[-\hbar\cdot \alpha_{i\bar{j}} d\bar{z}^j \otimes y^i,  \omega_{k \bar{l}} (\iota_{V}dz^k) \otimes \bar{y}^{l}  \right]_{\star}\\
		= &-\hbar\cdot \alpha_{i\bar{j}} dz^i (\iota_V d\bar{z}^j)  -\frac{2\pi}{\sqrt{-1}}\cdot\iota_V \omega + \hbar\cdot\alpha_{i\bar{j}} (\iota_V dz^i) d\bar{z}^j \\
		=&\frac{2\pi\hbar}{\sqrt{-1}}\cdot \iota_V \left(-\frac{1}{\hbar}\cdot\omega+\alpha\right)=\frac{2\pi\hbar}{\sqrt{-1}}\cdot\iota_VK.
		\end{align*}
	Therefore we see that we have to add the anti-derivative $\mu_V$ of $-\frac{2\pi\hbar}{\sqrt{-1}}\cdot\iota_VK$ to $\beta_V$ to achieve flatness locally, i.e. $D_{\alpha}(\mu_V+\beta_V) =d_X{\mu_V}+\frac{2\pi\hbar}{\sqrt{-1}}\cdot\iota_VK= 0$. Flatness can be achieved globally exactly when anti-derivative $\mu_V$ of $-\frac{2\pi\hbar}{\sqrt{-1}}\cdot\iota_VK$ can be solved globally. 
\end{proof}

\begin{prop}
	Suppose $\star$ is a Wick type star product on $X$ with Karabegov form $-\frac{1}{\hbar}\omega+\alpha$ which is admissible, i.e., $\alpha$ is a polynomial in $\hbar$. Suppose a vector field $V$ is a derivation with respect to $\star$ and has a quantum Hamiltonian. Then there exist quantum Hamiltonians which are quantizable functions on $X$. 
\end{prop}
\begin{proof}
	We only need to observe that the flat section associated to a quantum Hamiltonian is given by $\beta_V$ in equation \eqref{equation: flat-section-modulo-constant-term}. Since the term $\gamma$ is at most degree $1$ in $\overline{\W}_X$ and $\alpha$ is a polynomial in $\hbar$, it follows that $\beta_V$ is also a polynomial in $\hbar$ and $\overline{\W}_X$. Thus a quantum Hamiltonian must be quantizable.  
\end{proof}

%A most important example of Wick type star product is the Berezin-Toeplitz deformation quantization. 
By Proposition \ref{proposition: vector-field-derivation}, a vector field is a derivation of the Berezin-Toeplitz star product if and only if $V$ is compatible with both $J$ and $\omega$. 

\begin{corollary}\label{corollary: infinitesimal-symmetry-Berezin-Toeplitz}
	A vector field  $V$ on a prequantizable K\"ahler manifold $X$ acts as a derivation with respect to the Berezin-Toeplitz star product if and only if it is both symplectic and Killing. If an associated quantum Hamiltonian exists globally on $X$, then there exist quantum Hamiltonians which are degree $1$ quantizable functions. 
\end{corollary}

\subsection{Which degree $1$ quantizable functions arise from symmetry}
\

In this subsection we restrict ourself to Berezin-Toeplitz deformation quantization and letting $\alpha = Ric_X$. We have seen that quantum Hamiltonians determined by the equation$d\mu_V = \iota_V(\omega -\hbar Ric_X)$ gives rise to a large class of degree $1$ formal quantizable functions. In this subsection, we give an answer to the following converse question: which degree $1$ quantizable functions are quantum Hamiltonians?

\begin{prop}\label{proposition: Criterion-Hamiltonian-vector-field-preserve-complex-structure}
	For a real function $f$ on $X$, its associated symplectic vector field preserves the complex structure if and only if it satisfies the following condition:
	$$
	\nabla^{0,1}\left(\frac{\partial f}{\partial\bar{z}^j}\cdot\bar{y}^j\right)=0,
	$$
	which is also equivalent to $V_f^{1,0}$ being holomorphic. 
\end{prop}
\begin{proof}
	The above condition for $f$ can be written as
	\begin{equation}\label{equation: 0-2-term-flat-section-function}
		\begin{aligned}
			\nabla^{0,1}\left(\frac{\partial f}{\partial\bar{z}^j}\cdot\bar{y}^j\right)=&\frac{\partial^2 f}{\partial\bar{z}^j\partial\bar{z}^k}\cdot d\bar{z}^k\otimes\bar{y}^j+\frac{\partial f}{\partial\bar{z}^j}\cdot \nabla^{0,1}(\bar{y}^j)\\
			=&\frac{\partial^2 f}{\partial\bar{z}^j\partial\bar{z}^k}\cdot d\bar{z}^k\otimes \bar{y}^j+\frac{\partial f}{\partial\bar{z}^j}\cdot\left(\Gamma_{\bar{k}\bar{l}}^{\bar{j}}\cdot d\bar{z}^k\otimes\bar{y}^l\right)\\
			=&0. 
		\end{aligned}
	\end{equation}
	Recall that the Hamiltonian vector field of the function $f$ is given by the following explicit formula:
	$$
	V:=\frac{2\pi}{\sqrt{-1}}\cdot\omega^{i\bar{j}}\left(\frac{\partial f}{\partial z^i}\frac{\partial}{\partial\bar{z}^j}+\frac{\partial f}{\partial \bar{z}^j}\frac{\partial}{\partial z^i}\right).
	$$
	Since the function $f$ is real, to show that this vector field preserves the complex structure, we only need to show that $\mathcal{L}_V(dz^i)$ is still a form of type $(1,0)$. Using Cartan's magic formula, we have the following explicit computation:
	\begin{align*}
		\frac{\sqrt{-1}}{2\pi}\mathcal{L}_V(dz^i)=&\frac{\sqrt{-1}}{2\pi}\left(d(\iota_V\lrcorner dz^i)+\iota_V\lrcorner d(dz^i)\right)\\
		=&d\left(\omega^{i\bar{j}}\frac{\partial f}{\partial\bar{z}^j}\right)\\
		=&\partial\left(\omega^{i\bar{j}}\frac{\partial f}{\partial\bar{z}^j}\right)+\bar{\partial}\left(\omega^{i\bar{j}}\frac{\partial f}{\partial\bar{z}^j}\right)\\
		=&\partial\left(\omega^{i\bar{j}}\frac{\partial f}{\partial\bar{z}^j}\right)+\left(\omega^{i\bar{j}}\frac{\partial^2 f}{\partial\bar{z}^j\partial\bar{z}^k}d\bar{z}^k+\frac{\partial\omega^{i\bar{j}}}{\partial\bar{z}^k}d\bar{z}^k\cdot\frac{\partial f}{\partial\bar{z}^j}\right)\\
		=&\partial\left(\omega^{i\bar{j}}\frac{\partial f}{\partial\bar{z}^j}\right)+\left(\omega^{i\bar{j}}\frac{\partial^2 f}{\partial\bar{z}^j\partial\bar{z}^k}d\bar{z}^k+\omega^{i\bar{l}}\cdot\Gamma_{\bar{k}\bar{l}}^{\bar{j}}d\bar{z}^k\cdot\frac{\partial f}{\partial\bar{z}^j}\right).
	\end{align*}
	By comparing the rightmost term with equation \eqref{equation: 0-2-term-flat-section-function}, it is easy to see that $\mathcal{L}_V(dz^i)$ is still of type $(1,0)$ if and only if the condition for $f$ in the Lemma holds. 
	%	Thus the $V$ is indeed a Hamiltonian Killing vector field, which induces our desired one-to-one correspondence. 
\end{proof}

The main result in this subsection is the following:
\begin{theorem}\label{proposition: real-quantizable-function-Hamiltonian}
	A formal quantizable function $f$ on $X$ is a quantum Hamiltonian of a Hamiltonian Killing vector field on $X$ with respect to the Berezin-Toeplitz deformation quantization if and only if it is a real (up to a formal constant in $\C[[\hbar]]$) degree $1$ formal quantizable function. 
\end{theorem}
\begin{proof}
	Suppose $f$ is a real degree $1$ formal quantizable function, by Proposition \ref{proposition: degree-1-formal-quantizable-function}, it must be of the following form:
	$$
	f=f_0-\hbar\cdot (\frac{1}{4\pi}\Delta f_0 +c),
	$$
	where $f$ must be a classical Hamiltonian associated to $V$, i.e., $V=V_f$. Since $f$ is real, then so is $f_0$. By Proposition \ref{proposition: Criterion-Hamiltonian-vector-field-preserve-complex-structure}, $V_f$ preserves the complex structure on $X$. From the proof of Proposition \ref{proposition: degree-1-formal-quantizable-function}, we have 
	$$
	\frac{1}{4\pi}\frac{\partial\Delta f_0}{\partial\bar{z}^l} - \frac{\partial f_0}{\partial\bar{z}^j}\cdot\omega^{k\bar{j}}\cdot (Ric_X)_{k\bar{l}} = 0,
	$$
	or $\bar{\partial} (\frac{1}{4\pi}\Delta f_0)  = \iota_{V^{1,0}} Ric_X$. Taking real part of the equation we get $\iota_V Ric_X = d(\frac{1}{4\pi}\Delta f_0)$. 
	
	For the other direction, if $f$ is the quantum Hamiltonian of a vector field $V$, then we have $df =-\frac{2\pi\hbar}{\sqrt{-1}}\cdot\iota_VK_{BT}$ which forces it to be a degree $1$ formal quantizable function of the form $f=f_0-\hbar\cdot (\frac{1}{4\pi}\Delta f_0 +c)$ and the leading term $f_0$ is its classical Hamiltonian. Thus $f_0$ must be real (up to a constant) and so is $\Delta f_0$. 
\end{proof}

\begin{remark}
	This proposition shows that our definition of quantizable functions is much more general than the old notion of quantizable functions in geometric quantization (see for instance, section 4 in \cite{bordemann1991gl}). There a function is called quantizable if its associated Hamiltonian vector fields preserves the K\"ahler polarization, or equivalently the $(1,0)$-part of its Hamiltonian vector fields is holomorphic. (Actually, in this definition, it is implicitly required that the function is real since otherwise, only the requirement on $(1,0)$ part of the Hamiltonian vector field is not enough.) Thus Proposition \ref{proposition: real-quantizable-function-Hamiltonian} is saying that those functions are exactly degree $1$ real quantizable functions under the Definition in \cite{CLL23}.
\end{remark}

\begin{remark}
It is known in symplectic geometry, that symplectic vector fields corresponds to smooth functions (modulo constant functions). We can say that smooth functions corresponds to symmetries of the classical mechanical systems. 

Proposition \ref{proposition: real-quantizable-function-Hamiltonian} implies that we have a quantum analogue of this correspondence. From a physical point of view, a symplectic Killing vector field describes an infinitesimal symmetry of a K\"ahler manifold preserving both the algebraic structure of quantum observables, and the Hilbert spaces which depends on the complex polarization. These are the symmetries for the quantum mechanical system encoded in $X$. 
It is easy to see from Proposition \ref{proposition: real-quantizable-function-Hamiltonian} that there is a one-to-one correspondence between symplectic Killing vector fields and real degree $1$ quantizable functions (modulo constant functions) on a K\"ahler manifold $X$. 
\end{remark}

\section{Quantum symmetry in deformation and geometric quantization}\label{section: quantum-moment-map}

Suppose we have a Lie group $G$ (with Lie algebra $\mathfrak{g}$) acting on a manifold $X$ preserving certain geometric structures, then the infinitesimal symmetries are described in terms of Lie algebra homomorphisms from $\mathfrak{g}$ to vector fields on $X$. The infinitesmal $\mathfrak{g}$ symmetry can act as symmetry of corresponding quantum mechanical system, including the algebra of quantum observables (deformation quantization) and the Hilbert spaces (geometric quantization). 

We begin with the corresponding infinitesimal symmetries on quantum observables. This is  mathematically described as the notion of quantum moment map, introduced in \cite{Xu}.  Explicitly, suppose a Lie group $G$ acts on a symplectic manifold $X$ preserving a deformation quantization $\left(C^\infty(X)[[\hbar]],\star\right)$ in the sense that for any $\Phi\in G$, there is
$$
\Phi^*(f)\star\Phi^*(g)=\Phi^*(f\star g).
$$
Infinitesimally, this symmetry on deformation quantization is described as follows:
\begin{definition}\label{definition: quantum-moment-map}
	Let $X,G$ be as above. Suppose $(C^\infty(X)[[\hbar]],\star)$ is a $G$-invariant deformation quantization of $X$. Then a {\em quantum moment map} is a homomorphism of Lie algebras
	$$
	\mu_\hbar: \mathfrak{g}\rightarrow C^\infty(X)[[\hbar]],
	$$
	i.e., for $g,h\in\mathfrak{g}$, there is $\mu_\hbar([g,h])=\frac{1}{\hbar}[\mu_\hbar(g),\mu_\hbar(h)]_\star=\frac{1}{\hbar}\left(\mu_\hbar(g)\star\mu_\hbar(h)-\mu_\hbar(h)\star\mu_\hbar(g)\right)$. 	And that for every $g\in\mathfrak{g}$, we have the equality $\mathcal{L}_{V_g}=\frac{1}{\hbar}[\mu_\hbar(g),-]_\star$ for formal smooth functions $C^\infty(X)[[\hbar]]$; here $V_g$ denotes the vector field associated to $g \in \mathfrak{g}$. 
\end{definition}
\begin{remark}\label{remark: normalization-quantum-moment-map}
	%In the original definition of quantum moment maps in \cite{Xu}, the bracket on the Lie algebra $\mathfrak{g}$ is normalized with a factor $\hbar$. Here we revise this factor to $1/\hbar$ on $C^\infty(X)[[\hbar]]$. As we will see later, this will be more compatible with geometric quantization. 
	If we write a quantum moment map as
	$$
	\mu_\hbar=\mu_0+O(\hbar),
	$$
	then for any $\xi\in\mathfrak{g}$, the function $\frac{2\pi}{\sqrt{-1}}\cdot\mu_0(\xi)$ is a (classical) Hamiltonian associated to the vector field $V_\xi$. 
\end{remark}

 It is shown in \cite{Xu} that in general there are certain cohomological obstructions to the existence of quantum moment maps. For Wick type star product on (pseudo-)K\"ahler manifold, similar obstructions are also carefully studied in \cite{Muller-Neumaier}. 
 
 There are two sorts of cohomological obstructions to the existence of quantum moment maps, of which we give a quick review. First of all, it is easy to see that for any element  $g\in\mathfrak{g}$, its image under the quantum moment map is a quantum Hamiltonian associated to the vector field $V_g$. And we have seen from equation \eqref{equation: obstruction-quantum-hamiltonian} that there is a de Rham cohomology class which obstructs the existence of quantum Hamiltonian associated to $V_g$. Even if there are certain topological property of $X$ (for instance, $H^1_{dR}(X)=0$) which implies the vanishing of these obstructions, there are still the compatibilities between these quantum Hamiltonians. These compatibilities are described by a Lie algebra cohomology class in $H^2(\mathfrak{g},\C)[[\hbar]]$. To our best knowledge, little is known about this Lie algebra cohomology obstruction in general. 
 
 A possible scenario where these obstructions could vanish is the following: for every $g\in\mathfrak{g}$, there is a “canonical" quantum Hamiltonian associated to $g$, with the help of some additional structures. A natural choice is to consider the Hilbert space in geometric quantization on which the quantum observables could act naturally. To this end, we need a compatible symmetry, i.e., a $G$ action on the Hilbert spaces  compatible with the $G$ action on the quantum observables. In this circumstance, we can turn the infinitesimal action of any $g\in\mathfrak{g}$ canonically to a quantum observable without the ambiguity in $\C[[\hbar]]$, which implies the existence of a {\em quantum moment map}.

%\subsection{Quantum moment maps}
%\

%In the previous section, we know that the quantum Hamiltonian associated to a vector field is only unique up to a formal constant in $\C[[\hbar]]$. 

\subsection{Quantum symmetry on pre-quantizable K\"ahler manifolds}
\

In this subsection, we will only focus on the Berezin-Toeplitz quantization (both deformation and geometric quantization) on a pre-quantizable K\"ahler manifold $X$. In this geometric context, there are natural Hilbert spaces $\mathcal{H}_k:=H^0_{\bar{\partial}}(X, L^{\otimes k})$ of holomorphic sections of the (tensor power of) pre-quantum line bundle, on which the quantum observables act as operators. We will consider all (infinitesimal) symmetries on quantum observables and Hilbert spaces respectively.

We begin with the discussion of symmetries on the Hilbert spaces $\mathcal{H}_k:=H^0_{\bar{\partial}}(X, L^{\otimes k})$. Let $G$ be a group which acts on $X$ preserving both symplectic and complex structures. On one hand, this action can be lifted to the pre-quantum line bundle $L$ (and also its tensor powers) if and only if it preserves the symplectic form (i.e., the curvature of $L$). In this case, the group action $G$ preserves all the pre-Hilbert spaces $C^\infty(X, L^{\otimes k})$.  On the other hand, since there are enough holomorphic sections for $k>>0$, this $G$ action furthermore preserves all the Hilbert spaces if and only if $G$ also preserves the complex structure. 

Infinitesimally, a vector field acts as a symmetry on all these Hilbert spaces $\mathcal{H}_k$ if and only if it preserves both the symplectic and complex structures. We recall the explicit formula of the action of such a vector field $V_\xi$ on $H^0_{\bar{\partial}}(X, L^{\otimes k})$ for $\xi \in\mathfrak{g}$, which can be found in e.g. \cite{bordemann1991gl}. Let $s$ be a holomorphic section of $L^{\otimes k}$, define
\begin{equation}\label{equation: infinitesimal-symmetry-Hilbert-spaces}
	\beta_k(\xi):=\nabla_{V_\xi}^{\otimes k}-2\pi\sqrt{-1}k\cdot m_{\mu(\xi)}.
\end{equation}
Here $\nabla^{\otimes k}$ denotes the connection on the prequantum line bundle $L^{\otimes k}$, and $m_{\mu(\xi)}$ simply denotes the multiplication by $\mu(\xi)$.  Notice that this action is defined on the space of smooth sections of $L^{\otimes k}$ while it preserves the subspace of holomorphic sections. We denote the Lie algebra homomorphism by $\beta_k$:
\begin{equation}\label{equation: definition-beta-m}
	\beta_k: \mathfrak{g}\rightarrow gl(\mathcal{H}_k).
\end{equation}

 As to the symmetries on quantum observables (Berezin-Toeplitz deformation quantization), it is clear that if a vector field $V$ preserves a star product on $X$, it must be a derivation with respect to this star product. By Corollary \ref{corollary: infinitesimal-symmetry-Berezin-Toeplitz}, $V$ must be both symplectic and Killing. It follows that the corresponding group must act as automorphisms of $X$ as a K\"ahler manifold, which is identical to that of the symmetry group of the Hilbert spaces.

\subsection{Quantum moment map acts as symmetries for $\mathcal{H}_k$}
\

In this subsection, we will proceed to  show that for any $G$ action preserving the algebra of quantum observables and Hilbert space on a pre-quantizable $X$, the corresponding infinitesimal Lie algebra actions are compatible through the Bargmann-Fock action. 

More precisely, for every $\xi\in\mathfrak{g}$, we will define an associated smooth formal function as follows. First of all, equation \eqref{equation: infinitesimal-symmetry-Hilbert-spaces} gives rise to a canonical classical moment map $\mu(\xi)$. Secondly, we define a map 
$$
\mu_\hbar:\mathfrak{g}\rightarrow C^\infty(X)[[\hbar]]
$$
by the following equation 
\begin{equation}\label{equation: canonical-quantum-Hamiltonian}
	\mu_\hbar(\xi):=\mu(\xi)-\frac{\hbar}{4\pi} \Delta(\mu(\xi)).
\end{equation}
\begin{lemma}
For every $\xi\in\mathfrak{g}$, the formal function $\mu_\hbar(\xi)$ is a degree $1$ formal quantizable function.
\end{lemma}
\begin{proof}
	It is easy to see that if the vector field $V_\xi$ preserves the complex structure on $X$, the image under the classical moment map $\mu(\xi)$ must satisfy equation \eqref{equation: f-0}, and Proposition \ref{proposition: degree-1-formal-quantizable-function} implies that $\mu_\hbar$ is a degree $1$ formal quantizable function. 
\end{proof}

We can thus define the following degree $1$ level $k$ quantizable functions by taking the evaluation $\hbar=1/k$:
\begin{equation}
	\mu_k: \mathfrak{g} \rightarrow C^\infty_{\alpha,k}(X). 
\end{equation}
\begin{remark}
	We will show later that the map $\mu_\hbar$ is a {\em quantum moment map} associated to this symmetry on quantum observables. And we will call $\mu_k$'s the level $k$ quantum moment map and justify it later. 
\end{remark}
As we have seen in section \ref{section: Hilbert-spaces-quantizable-functions}, a quantizable function of level $k$ acts on $H^0_{\bar{\partial}}(X, L^{\otimes k})$ via the Bargmann-Fock action. And we have the first main result:
%and call it the {\em quantum moment map of level $k$}. 

%We should see how this infinitesimal symmetry relate to symmetry of geometric quantization. 

%In particular, equation \eqref{equation: infinitesimal-symmetry-Hilbert-spaces} implies that for any vector field $V\in\mathfrak{g}_X$, there exists a canonical classical moment map $\mu(V)$ associated to $V$.  Since $V$ preserves the complex structure, the classical Hamiltonian $\mu(V)$ satisfies the condition in Proposition \ref{proposition: degree-1-formal-quantizable-function}, thus the following function in equation \eqref{equation: canonical-quantum-Hamiltonian} is a degree $1$ formal quantizable function. 

%In the previous sections, we see that the quantum moment map $\mu_k : \mathfrak{g} \rightarrow \mathcal{C}^{\infty}_{\alpha,k}(X)$ can be viewed as an infinitesmal symmetry for deformation quantization. On the other hand, we have $\beta_k$ in equation \eqref{equation: definition-beta-m} which serves as infinitesmal symmetry for geometric quantization. 

\begin{theorem}\label{theorem: compatibility-quantum-moment-map-differetial-operator}
	For each level $k>0$, there exists a map $\mu_k:\mathfrak{g}\rightarrow C_{\alpha,k}^\infty(X)$, such that the following diagram commutes:
	\begin{equation}\label{equation: commutative-diagram-differential-operator}
		\xymatrix{ & C^\infty_{\alpha,k}(X) \ar[dr]^{\mathcal{D}}  & \\
			\mathfrak{g} \ar[ur]^{{k\cdot\mu}_k}\hspace{2mm}  \ar[rr]_{\beta_k} & & gl(\mathcal{H}_k)
		}
	\end{equation}
 Each $\mu_k$ is the evaluation of a quantum moment map $\mu_\hbar$ at $\hbar=1/k$. Equivalently,  (infinitesimal) symmetries on the Hilbert spaces $H^0_{\bar{\partial}}(X, L^{\otimes k})$ can be equivalently described via the Bargmann-Fock action of the level $k$ quantum moment map $\mu_k$.
\end{theorem}
\begin{proof}
	Given any $\xi\in\mathfrak{g}$, we let $V_\xi$ denote the associated vector field on $X$. Let $\mu(\xi)$ denote the function defined by equation \eqref{equation: infinitesimal-symmetry-Hilbert-spaces} (classical moment map). 
	Next we define a formal function $\mu_\hbar(\xi)$ by
	$$
	\mu_\hbar(\xi)=\frac{2\pi}{\sqrt{-1}}\left(\mu-\frac{1}{4\pi}\hbar\Delta\mu\right). 
	$$
	
	 Notice that the vector field associated to $\mu(\xi)$ is then explicitly 
	$$
	V_{\mu(\xi)}=\frac{2\pi}{\sqrt{-1}}\left(\frac{\partial\mu}{\partial z^i}\omega^{i\bar{j}}\frac{\partial}{\partial\bar{z}^j}+\frac{\partial \mu}{\partial\bar{z}^j}\omega^{\bar{j}i}\frac{\partial}{\partial z^i}\right).
	$$
	By taking the evaluation $\hbar=1/k$, we obtain the function $\mu_k(\xi)$. Since $V_\xi=V_{\mu_\xi}$ preserves the complex structure on $X$, it follows from Proposition \ref{proposition: degree-1-formal-quantizable-function} that $\mu_\hbar(\xi)$ is a degree $1$ formal quantizable function. 
	
	%It is clear that the $\mathfrak{g}$ action on the Hilbert spaces $\mathcal{H}_k$ is via a first order differential operator.  We can look at the flat sections associated to both the quantum moment map and the holomorphic section in $\mathcal{H}_k$. 
	
	 Let $s$ be a  holomorphic section of $L^{\otimes k}$, whose associated flat section is locally of the form $O_s=O_g\cdot e^{k\cdot\beta}\otimes e_{L^{k}}$ as in equation \eqref{equation: flat-section-Hilbert-spaces}. The horizontal arrow in diagram \eqref{equation: commutative-diagram-differential-operator} is 
	 \begin{align*}
	 	\beta_k(\xi)(s)=&\nabla_{V_\xi}^{\otimes k}(s)-2\pi\sqrt{-1}k\cdot\mu(\xi)\cdot s\\
	 	=&V_\xi(g)\cdot s+k\cdot\frac{\partial\mu}{\partial\bar{z}^j}\omega^{\bar{j}i}\frac{\partial\rho}{\partial z^i}\cdot s-2\pi\sqrt{-1}k\cdot\mu(\xi)\cdot s. 
	 \end{align*}

	 On the other hand, the composition of the upper arrows in diagram \eqref{equation: commutative-diagram-differential-operator} is 
	 $$
	 \mathcal{D}\circ k\cdot\mu_k(\xi)=\sigma\left(O_{k\cdot\mu_k(\xi)}\circledast_k O_{s}\right).
	 $$
	  %To find the symbol of $O_{\mu_k(\xi)}\circledast_k O_{s}$, we need to look at the flat section corresponding to  ${\mu_k(\xi)}$ in more details. 
	  %The components of $O_{\mu_k(\xi)}$ in the anti-holomorphic Weyl bundle $\overline{\mathcal{W}}_X$ is given by
	%$$\theta = f_0-\frac{1}{2k} \Delta(f_0) + \frac{\partial f_0}{\partial \bar{z}^j} \bar{y}^j,$$
%	in local coordinates, which is at most degree $1$ in $\bar{y}^j$'s. From \cite[Proposition 2.24.]{CLL23}, we have $O_{\mu_k(\xi)} = \sum_{m\geq 0 } \left( \tilde{\nabla}^{1,0} \right)^m (\theta)$. 
		%According to the discussion in Remark \ref{remark: normalization-quantum-moment-map},  vector field $V_{f_0}$ is given by 
%$$V_{f_0}=\frac{2\pi}{\sqrt{-1}}\left(\frac{\partial f_0}{\partial z^i}\omega^{i\bar{j}}\frac{\partial}{\partial\bar{z}^j}+\frac{\partial f_0}{\partial\bar{z}^j}\omega^{\bar{j}i}\frac{\partial}{\partial z^i}\right)$$
%	It follows that the corresponding covariant derivative on $e_L^k$ is explicitly
%	\begin{align*}
%		\nabla^{\otimes k}_{V_{f_0}}(e_L^k)=&\\
%		=&
%	\end{align*}
	Since we know $O_{k\cdot\mu_k(\xi)}$ has at most polynomial degree $1$ terms in the anti-holomorphic Weyl bundle $\overline{\W}_X$, the only terms in the expression of $O_{k\cdot\mu_k(\xi)} $ that could contribute to the the symbol of $O_{k\cdot\mu_k(\xi)}\circledast_k O_{s}$ are
	$$
	-2\pi\cdot k\sqrt{-1}\left(\mu(\xi)-\frac{1}{4\pi k} \Delta(\mu(\xi)) + \frac{\partial\mu(\xi)}{\partial \bar{z}^j} \bar{y}^j + \frac{\partial^2\mu(\xi)}{\partial z^i \partial \bar{z}^j} y^i \bar{y}^j\right), 
	$$
	\begin{comment}
	Since we know $O_{\mu_k(\xi)}$ has at most degree $1$ in the anti-holomorphic Weyl bundle, We now look at the following terms in $O_{\mu_k(\xi)}$:
	$$
	\sum_{m\geq 0}\left(\tilde{\nabla}^{1,0}\right)^m\left(\frac{\partial^2 f_0}{\partial z^i\partial\bar{z}^j}\cdot y^i\bar{y}^j\right)
	$$
	We have the following equality:
	$$
	\tilde{\nabla}^{1,0}\left(\frac{\partial^2 f_0}{\partial z^i\partial\bar{z}^j}\cdot y^i\bar{y}^j\right)=R_{i\bar{j}k\bar{l}}\cdot V_{\xi}(d\bar{z}^j)\cdot y^iy^k\bar{y}^l. 
	$$
	For the above term, there is 
	$$
	\omega^{k\bar{l}}\cdot R_{i\bar{j}k\bar{l}}\cdot V_{\xi}(d\bar{z}^j)\cdot y^i=(\text{Ric}_X)_{i\bar{j}}\cdot V_{\xi}(d\bar{z}^j).
	$$
	This implies that these two types of terms in $O_{\mu_k(\xi)}\circledast_k O_{s}$ will cancel with each other, 
	\end{comment}
	and we have 
	\begin{align*}
		&\sigma(O_{k\cdot\mu_k(\xi)}\circledast_k O_{s})\\
		=&-2\pi\cdot k\sqrt{-1} \left(\mu(\xi)\cdot s - \frac{1}{4\pi k}(\Delta\mu(\xi)) \cdot s+ \sigma\left( ( \frac{\partial\mu(\xi)}{\partial\bar{z}^j}\bar{y}^j +  \frac{\partial^2\mu(\xi)}{\partial z^i \partial \bar{z}^j} y^i \bar{y}^j)\circledast O_s \right)\right)\\
		=&-2\pi\cdot k\sqrt{-1}\left(\mu(\xi)-\frac{1}{4\pi k}(\Delta\mu(\xi)) \right)\cdot s\\
		&+\sigma\left( -2\pi\sqrt{-1}\big(\frac{\partial\mu(\xi)}{\partial\bar{z}^j} \omega^{\bar{j}k}\frac{\partial}{\partial y^k} +\frac{\partial^2\mu(\xi)}{\partial z^i \partial \bar{z}^j} \omega^{\bar{j}k} \frac{\partial}{\partial y^k} \circ m_{y^i} \big) (O_g\cdot e^{k\cdot\beta} )\otimes O_{e_{L^k}}\right)\\
		=&-2\pi\cdot k\sqrt{-1}\left(\mu(\xi)-\frac{1}{4\pi k}(\Delta\mu(\xi)) \right)\cdot s\\
		&-2\pi\sqrt{-1}\left(\omega^{\bar{j}k}\frac{\partial\mu(\xi)}{\partial\bar{z}^j}\left(\frac{\partial g}{\partial z^k} +g\cdot k\frac{\partial \rho}{\partial z^k}\right)\otimes e_{L^k} +\frac{1}{4\pi}\Delta\mu(\xi)\cdot s\right)\\
		=&-2\pi\cdot k\sqrt{-1}\cdot\mu(\xi)\cdot s +V_\xi(g)\cdot e_{L^k}+ g\cdot\nabla^{\otimes k}_{V_\xi}(e_{L^k})\\
		=&-2\pi\cdot k\sqrt{-1}\cdot\mu(\xi)\cdot s+\nabla^{\otimes k}_{V_\xi}(s) = \beta_k(\xi),
	\end{align*}
	which says that the Bargmann-Fock action of $k\cdot\mu_k(\xi)$ on $s$ exactly coincides with the corresponding infinitesimal symmetry. Here we have used the explicit formula of $O_{e_{L^k}}$ in equation \eqref{equation: flat-section-Hilbert-spaces}.  
	
	Next we show that $\mu_\hbar$ is a quantum moment map for the Berezin-Toeplitz deformation quantization. The commutative diagram \eqref{equation: commutative-diagram-differential-operator} in Theorem \ref{theorem: compatibility-quantum-moment-map-differetial-operator} implies that for all $\xi_1,\xi_2\in\mathfrak{g}$, there is
	$$
	k\cdot[\mu_k(\xi_1),\mu_k(\xi_2)]_{\star_k}=\mu_k([\xi_1,\xi_2]).
	$$
	(Here we are also using the fact that $\beta_k$ is a Lie algebra homomorphism).
	 Since this holds for all level $k>0$, we can turn these $k$'s to the formal variable and $\mu_\hbar: \mathfrak{g}\rightarrow C^\infty(X)[[\hbar]]$ is a quantum moment map. 
\end{proof}

\begin{corollary}
	Let $G$ be a Lie group which acts on a pre-quantizable K\"ahler manifold $X$ preserving the K\"ahler and complex structures, then there exists a quantum moment map for the Berezin-Toeplitz deformation quantization. 
\end{corollary}

%\begin{remark}
%	We can check that the classical part of $\mu_\hbar$, i.e., $\mu$ is indeed a classical moment map. 
%\end{remark}

%We begin with the existence of quantum moment maps for a Hamiltonian Lie Group action $G$ on $X$ with classical moment map $\mu : X \rightarrow \mathfrak{g}^*$, using our understanding of quantum Hamiltonian functions. 

%Theorem \ref{theorem: compatibility-quantum-moment-map-differetial-operator} and Proposition \ref{proposition: canonical-quantum-Hamiltonian-quantum-moment-map} tells us although quantum moment map is defined in deformation quantization, their images also acts on the Hilbert spaces in geometric quantization as differential operators in a compatible way. 

%Since the $G$ action is assumed to be Hamiltonian, the vector field $V$ preserves the pre-quantum line bundle on $X$ (****$H^1(\mathcal{O}_X)=0?$****). Furthermore, $V$ can act infinitesimally on the Hilbert spaces in geometric quantization $H^0_{\bar{\partial}}(X, L^{\otimes k})$ by the following explicit form: for any $\xi\in\mathfrak{g}_X$, we let $V_\xi$ denote the associated vector field on $X$. The associated $\mathfrak{g}$ action on the Hilbert spaces is defined as follows: 

\subsection{Symmetry in Berezin-Toeplitz geometric quantization}\label{section: Berezin-Toeplitz-quantization}
\

Another more widely used quantization method is the so called Berezin-Toeplitz quantization, in which every smooth function is quantized to a Berezin-Toeplitz operator acting on the Hilbert spaces $\mathcal{H}_k$. 
We first briefly recall the {\em Berezin-Toeplitz quantization}. For a level $k>0$, we let $P_k: \Gamma(X, L^{\otimes k})\rightarrow\mathcal{H}_k$ denote the orthogonal projection from smooth sections of $L^{\otimes k}$ to the subspace of holomorphic sections. For every smooth function $f\in C^\infty(X)$, we define the Berezin-Toeplitz operator on $\mathcal{H}_k$ as
$$
T_{f,k}:=P_k\circ m_f. 
$$
Here $m_f$ denotes the multiplication by the function $f$. The asymptotic formula as $k\rightarrow \infty$ of the composition of two Berezin-Toeplitz operators induces the same {\em Berezin-Toeplitz star product} as our earlier construction described in section \ref{section: Fedosov-quantization} using the Karabegov form in equation \eqref{equation: Karabegov-form-Berezin-Toeplitz}. 

We first describe how these vector fields corresponding to infinitesimal symmetries acts on the analytically defined Berezin-Toeplitz operators. 
\begin{theorem}
	The $G$-action on the K\"ahler manifold $X$ naturally induces a quantum $G$-action on Berezin-Toeplitz operators. Explicitly, let $\xi\in\mathfrak{g}$ be an infinitesimal action on $X$ (and also on the spaces $\mathcal{H}_k$ of holomorphic sections). There is the following equality:
	$$
	[\beta_k(\xi), T_f]=T_{[k\cdot\mu_k(\xi),f]_{\star_k}}.
	$$
\end{theorem}
\begin{proof}
	Recall that the Berezin-Toeplitz operator $T_f=P_k\circ m_f\circ P_k$, where $P_k$ denotes the orthogonal projection and $m_f$ is the multiplication by the function $f$. Since the Hermitian metric on $L$ and the volume form on $X$ are both $G$-invariant, the orthogonal projections $P_k$ must be $G$-equivariant, which implies the commutativity $[\beta_k(\xi),P_k]=0$, for any $\xi\in\mathfrak{g}$. Thus we have 
	\begin{align*}
		[\beta_k(\xi), T_f]=&P_k\circ [\nabla^{\otimes k}_{V_\xi}- 2\pi\sqrt{-1}m_{\mu(\xi)}, m_f]\circ P_k\\
		=&P_k\circ [\nabla^{\otimes k}_{V_\xi},m_f]\circ P_k=P_k\circ V_\xi(f)\circ P_k\\
		=& T_{V_\xi(f)}=T_{[k\cdot\mu_k(\xi),f]_{\star_k}}.
	\end{align*}
	Here we have used the definition of the quantum moment map in the last equality. 
\end{proof}

An interesting question is that if the differential operators in diagram \eqref{equation: commutative-diagram-differential-operator} is replaced by the differential operators by Berezin-Toeplitz operators, does the commutativity still hold? 
\begin{equation}\label{equation: commutative-diagram-Toeplitz-operators}
	\xymatrix
	{ & C^\infty(X) \ar[dr]^{T}  & \\
		\mathfrak{g} \ar[ur]^{{k\cdot\mu}_k}\hspace{2mm}  \ar[rr]_{\beta_k} & & gl(\mathcal{H}_k)
	}
\end{equation}

Here we give the following example which shows that the diagram \eqref{equation: commutative-diagram-differential-operator} is not commutative in general for non-compact K\"ahler manifolds.  
\begin{example}
	We consider the complex plane $\C$ with the standard complex and symplectic structure, with the $G=\mathbb{R}$ symmetry given by the translation along the $x$-axis. We now let $X_R$ be an open strip $\{z=x+\sqrt{-1}\cdot y: -R<y<R\}$ which still allows the $\mathbb{R}$-symmetry. Thus the quantum moment map for all these $X_R$ and $\C$ are the same. On all these K\"ahler manifolds, the prequantum line bundle is trivial with the hermitian inner product given by 
	$$
	\langle f,g\rangle=f\cdot\bar{g}\cdot e^{-k\cdot |z|^2}.
	$$
	A straightforward computation shows that even for the trivial holomorphic section $1$, the orthogonal projection for different $X_R$'s are obviously different, and it follows that the above diagram cannot be commutative for all radius $R$.
\end{example}

However, for compact K\"ahler manifolds the diagram \eqref{equation: commutative-diagram-Toeplitz-operators} is indeed commutative.  To prove this, we need the following lemma:
\begin{lemma}[Tuynman \cite{tuynman1987quantization}, see also \cite{bordemann1991gl}*{Proposition 4.1.}]\label{lem:tuynman_formula}
	Let $(X,\omega,J,L)$ be compact K\"ahler manifold with a pre-quantum line bundle, then we have 
	$$
	\Pi_k \circ \nabla^{\otimes k}_{V_f} =\frac{\sqrt{-1}}{2}\Pi_k \circ m_{\Delta(f)}
	$$
	acting on holomorphic sections $H^0_{\bar{\partial}}(X,L^{\otimes k})$, where $\Delta$ is the Laplacian with respect to the K\"ahler metric.
\end{lemma}We now state and prove the following main theorem of this section:
\begin{theorem}
	Let $X$ be a compact prequantizable K\"ahler manifold, and let $G$ act on $(X,L)$ as before. Then the diagram \eqref{equation: commutative-diagram-Toeplitz-operators} is commutative. 
\end{theorem}
\begin{proof}
	For any $\xi\in\mathfrak{g}$, we have 
	\begin{align*}
		T_{k\cdot{\mu}_k(\xi)}=&\Pi_k\circ(-2\pi k\cdot\sqrt{-1})\cdot (m_{\mu(\xi)}-\frac{1}{2}\cdot\frac{1}{2\pi k}\cdot m_{\Delta\mu(\xi)})\\
		=&\Pi_k\circ(-2\pi k\cdot\sqrt{-1}\cdot m_{\mu(\xi)}+\frac{\sqrt{-1}}{2}\cdot m_{\Delta\mu(\xi)}) \\
		=&\Pi_k\circ(-2\pi k\cdot\sqrt{-1}\cdot m_{\mu(\xi)}+\nabla^{\otimes k}_{V_{\xi}}) \\
		=&-2\pi k\cdot\sqrt{-1}\cdot m_{\mu(\xi)}+\nabla^{\otimes k}_{V_{\xi}}\\
		=&\beta_k(\xi),
	\end{align*} 
	as action on holomorphic sections $H^0_{\bar{\partial}}(X,L^{\otimes k})$. In the second equality, we have used Lemma \ref{lem:tuynman_formula}. 	
\end{proof}

\end{document}